# Noble-Abel / First-order virial equations of state for gas mixtures resulting of multiple condensed reactive materials combustion

Loann Neron [1,2] and Richard Saurel [1,2]

[1] RS2N SAS, 371 chemin de Gaumin, 83640 Saint-Zacharie, France

[2] Aix Marseille Univ, CNRS, Centrale Marseille, LMA UMR 7031, Marseille, France

## Abstract

The Noble-Abel (NA) equation of state (EOS) is widely used in interior ballistics of guns as well as rocket propulsion computations. Its simplicity and accuracy are key points for intensive computations with hyperbolic two-phase flow models considered in interior ballistics codes. An alternative is examined in the present contribution through a first-order virial (VO1) equation of state. Appropriate methods for the determination of related parameters, such as specific gas constant, covolume and condensed material energy for both formulations (NA and VO1) are presented. Combination of closed bomb vessel experiments and thermochemical code computations are needed. An extended VO1 EOS with temperature dependent specific heat is examined. Then extension to multiple reactive materials is addressed. Examples are examined for each formulation (NA and VO1) and comparisons are done with the Becker-Kistiakowsky-Wilson (BKW) EOS as reference. Several conclusions emerged. First, consideration of specific heat temperature dependance in interior ballistics of guns computations appeared insignificant. Second, VO1 appeared more accurate than NA, particularly when gas density comes out of the range used for the EOS parameters determination. Last, regarding mixtures of condensed reactive materials, producing burnt gas mixtures, NA appeared again less accurate than VO1. However, its formulation is explicit, while VO1 requires numerical solving of a non-linear equation, with consequences on computational cost.

**Keywords :** Equation of state – Compressible flows – Dense gases – Ballistics of guns

(1) Corresponding author: loann.neron@etu.univ-amu.fr, loann.neron@rs2n.eu

(2) richard.saurel@univ-amu.fr

## 1. Introduction

The present contribution deals with the examination of various reduced EOS for interior ballistics of guns or solid rocket engines computations. In these areas, single phase and multiphase compressible flow models are needed, such as for example Baer and Nunziato (1986), Saurel and Abgrall (1999) and Saurel et al. (2017) models. Hyperbolic flow models only are mentioned, as wave propagation has importance in these applications. In these models, the thermodynamic closure for the solid phase may be addressed through stiffened gas (Le Métayer et al., 2004) or Noble-Abel stiffened gas (Le Métayer and Saurel, 2016) EOS for example. More sophisticated formulations are possible options, such as Cochran and Chan (1979) EOS. Regarding the gas phase, sophisticated formulation such as Becker-Kistiakowsky-Wilson (BKW) (Becker, 1921), (Kistiakowsky and Wilson, 1941) or high-order virial expansions (Heuzé, 1986) may be used as well. However, when handling multidimensional computations with non-equilibrium hyperbolic multiphase flow models, at least three key features must be present in the thermodynamic closure:

- The thermodynamic model must be obviously accurate.
- The EOS must be convex, otherwise hyperbolicity may be lost locally, resulting in numerical difficulties and unphysical wave propagation.
- The formulation must be computationally efficient. BKW as well as high-order virial EOS involve much higher cost, when coupled to hydrodynamic codes compared to simpler formulations for two main reasons. The first reason lies with non-explicit transformations, inherent to these EOS, to extract state variables, for instance temperature computation from internal energy. The second reason is specific to hyperbolic multiphase models where stiff pressure relaxation solvers are used (Lallemand et al., 2005). The latter, in conjunction with complicated EOS, involves a non-linear system of equations to be solved at each time step resulting in extra computational cost.

The present contribution follows these lines. Reduced EOS, such as Noble-Abel (NA) and first-order virial (VO1) are considered, thanks to their simplicity. Their accuracy is examined carefully versus more sophisticated formulations in closed bomb situations.

Reduced and sophisticated EOS for the gas phase have already been considered in many contributions, such as for example Farrar and Leeming (1983) or Wang et al. (1985). Most interior ballistics codes use the NA EOS (Farrar and Leeming, 1983), as it is simple, computationally efficient and convex. However, accuracy of this EOS is restricted to the density range used for the determination of its various parameters. When the EOS is used out of this density range, significant deviations appear as will be shown in this contribution.

In advanced applications, it is difficult and unsafe to perform closed bomb experiments at density levels representative of guns operative conditions. However, modern and future gun systems are expected to reach very high gas-density levels, exceeding 500 $kg/m^3$. When NA EOS is used with parameters determined in the range 100-200 $kg/m^3$, computed pressure at density 500 $kg/m^3$ becomes inaccurate. Therefore, peak pressure uncertainties are present and performance optimization of these systems becomes even more difficult.

The same goal of high-performance guns with very high energetic material loading rate is investigated in van Driel et al. (2017), Daniel et al. (2017), Yang et al. (2020), Conroy et al. (2021) under various aspects, ranging from material processing to computational geometry.

The issue related to appropriate design of the propellant charge is investigated in Conroy et al. (2021). A specific method to design innovative gun loadings has been done at RS2N, France (www.rs2n.eu). However, issues related to thermodynamics of high-density gases seems omitted or considered through expensive closed bomb vessels able to operate at very high pressures.

In the present work, it is shown that first-order virial (VO1) equation of state is a better alternative than NA EOS when the gas density is extrapolated outside the range used for thermodynamic parameters determination, while remaining simple, computationally fast and convex.

Temperature dependance of the specific heat in VO1 EOS is examined with the aim of improving its accuracy. It is shown that in high pressure conditions, such sophistication is not relevant.

Mixtures of reactive solid materials are then studied, as the main goal of the present contribution. Only mixtures of condensed energetic materials with same oxygen balance sign are considered. Their combustion results in gas mixtures for which reduced mixture thermodynamics is needed. Having the same oxygen balance sign, post combustion effects between the various gases are absent and therefore omitted in the present mixing rule. Both NA and VO1 EOS are considered for such extension. Mixture rules based on mixture internal energy and mixture specific volume definitions are used, with the assumption of well mixed gases, evolving in temperature and pressure equilibrium. Applying this method results in MNA and MVO1 EOS (M stands for mixture). Their accuracy is assessed against CHEETAH 2.0 computations (Fried and Souers, 1994), where the sophisticated BKW EOS is used as reference for mixtures of materials.

The paper is organized as follows. The formulations of the reduced EOS (NA and VO1) are described in Section 2. The methodology used to determine parameters of these EOS is presented in Section 3. Both EOS are tested on various examples in Section 4. Extension to temperature dependent specific heat in VO1 EOS is examined in Section 5. Extension to mixtures of reactive solid materials is addressed in Section 6. Various examples are studied in Section 7. Conclusions are given in Section 8.

## 2. Formulation of NA and VO1 EOS

The two reduced EOS considered in the present work are presented hereafter. Their compatibility with Maxwell's rules is examined as well.

### *2.1. The Noble-Abel equation of state*

The Noble-Abel thermal and caloric EOS read,

$$P(v,T)=\frac{RT}{v-b}, \tag{2.1}$$

$$e(T)=C_v(T-T_0)+e_0, \tag{2.2}$$

where P, T, v and e represent respectively the pressure, the temperature, the specific volume and the specific internal energy of the gas produced by the considered reactive material. $e_0$ is the reference energy of the gas and $T_0$ the reference temperature. R, b and $C_v$ are parameters characteristic of the thermodynamic properties of the material considered, representing respectively the specific gas constant, the covolume and the specific heat at constant volume. The specific gas constant is defined as the universal gas constant divided by the molar mass of the gas, $R=\hat{R}/\hat{W}$.

The following energy constant q is defined to simplify EOS expressions,

$$q=e_0-C_vT_0. \tag{2.3}$$

Introducing (2.3) in (2.1) and (2.2), NA thermal and caloric EOS become,

$$P(v,e)=\frac{R(e-q)}{C_v(v-b)}, \tag{2.4}$$

$$e(T)=C_vT+q. \tag{2.5}$$

Thermal and caloric EOS must fulfill a compatibility condition derived from Maxwell rules (Callen and Kestin, 1960). This condition expresses thermodynamic compatibility of (2.1) and (2.2). The first Maxwell rule reads,

$$\left(\frac{\partial}{\partial v}\left(\frac{\partial f}{\partial T}\right)_v\right)_T = \left(\frac{\partial}{\partial T}\left(\frac{\partial f}{\partial v}\right)_T\right)_v , \tag{2.6}$$

where $f = e - Ts$ is the Helmholtz free energy. With the help of Gibbs identity,

$$de = Tds - Pdv, \tag{2.7}$$

where s is the specific entropy, the following partial derivatives arise,

$$\begin{aligned} \left(\frac{\partial f}{\partial T}\right)_v &= -s, \\ \left(\frac{\partial f}{\partial v}\right)_T &= -P. \end{aligned} \tag{2.8}$$

By introducing (2.8) into the first Maxwell rule (2.6) and with the help of (2.7), the compatibility condition becomes,

$$\left(\frac{\partial e}{\partial v}\right)_T = T\left(\frac{\partial P}{\partial T}\right)_v - P. \tag{2.9}$$

From (2.1) and (2.2) the following partial derivatives are derived,

$$\begin{aligned} \left(\frac{\partial e}{\partial v}\right)_T &= 0, \\ \left(\frac{\partial P}{\partial T}\right)_v &= \frac{R}{v-b}. \end{aligned} \tag{2.10}$$

Inserting (2.10) in (2.9) results in,

$$\left(\frac{\partial e}{\partial v}\right)_T = \frac{RT}{v-b} - P = 0. \tag{2.11}$$

The compatibility condition (2.6) is consequently fulfilled, meaning that (2.1) and (2.2) form a compatible set of thermal and caloric EOS

The specific enthalpy is another useful function and needs to be formulated. From equation (2.1) the specific volume reads,

$$v(P,T) = \frac{RT}{P} + b. \tag{2.12}$$

The specific enthalpy is readily obtained as,

$$h(P,T) = e(P,T) + Pv(P,T) = (R + C_v)T + bP + q. \tag{2.13}$$

The specific heat at constant pressure is deduced from the specific enthalpy as,

$$C_p = \left(\frac{\partial h}{\partial T}\right)_P = R + C_v. \tag{2.14}$$

Relation (2.14) corresponds to Mayer's relation for NA EOS. Introducing the specific heat ratio $\gamma = C_p/C_v$,

$$\gamma - 1 = \frac{R}{C_v}, \tag{2.15}$$

the same relationship for ideal gases is obtained.

In compressible flow models, knowledge of the speed of sound is fundamental. It is defined as,

$$c^2 = -v^2 \left(\frac{\partial P}{\partial v}\right)_s = \frac{v^2 \left[P + \left(\frac{\partial e}{\partial v}\right)_P\right]}{\left(\frac{\partial e}{\partial P}\right)_v}. \tag{2.16}$$

For the NA EOS, it reads,

$$c^2(P, \rho) = \frac{c_{IG}^2}{1 - \rho b}, \tag{2.17}$$

where $c_{IG}^2 = \gamma P/\rho$ represents the square sound speed for an ideal gas, with the gas density $\rho = 1/v$. Expression of the specific entropy is given in Appendix A. Convexity of the NA equation of state is examined in Appendix B. It is based on the second, third and fourth Maxwell rules as well as other requirements.

### *2.2. The first-order virial equation of state*

The virial expansion (Onnes, 1902) reads,

$$Z = \frac{P}{\rho RT} = 1 + a\rho + b\rho^2 + c\rho^3 + ..., \tag{2.18}$$

where Z represents the compressibility factor and ρ the gas density. a, b, c, … are the virial coefficients, usually temperature dependent. In this work, only a first-order virial expansion with a constant virial coefficient is considered. The VO1 thermal and caloric equations of state are then,

$$P(\rho, T) = \rho RT(1 + a\rho), \tag{2.19}$$

$$e(T) = C_v T + q. \tag{2.20}$$

With the energy constant q defined as,

$$q = e_0 - C_v T_0. \tag{2.21}$$

It is convenient to express the pressure as a function of density and internal energy, similar to equation (2.4) for NA.

$$P(\rho, e) = \frac{\rho R}{C_v}(e - q)(1 + a\rho). \tag{2.22}$$

Again, thermal and caloric EOS must fulfill the compatibility condition (2.9), expressed as,

$$\left(\frac{\partial e}{\partial \rho}\right)_T = -\frac{T}{\rho^2}\left(\frac{\partial P}{\partial T}\right)_\rho + \frac{P}{\rho^2}. \tag{2.23}$$

The following partial derivatives are deduced from (2.19) and (2.20),

$$\left(\frac{\partial e}{\partial \rho}\right)_T = 0,$$
$$\left(\frac{\partial P}{\partial T}\right)_\rho = \rho R\left(1+a\rho\right). \tag{2.24}$$

Replacing the partial derivatives in (2.23) with (2.24), the compatibility condition appears satisfied, meaning that (2.19) and (2.20) form a compatible set of thermal and caloric EOS.

To determine the specific enthalpy, the density expressed as a function of pressure and temperature is needed. Equation (2.19) is a quadratic function with respect to density whose positive root is,

$$\rho\left(P,T\right) = \frac{-1+\sqrt{1+\frac{4aP}{RT}}}{2a}. \tag{2.25}$$

Expression of the specific enthalpy is then,

$$h\left(P,T\right) = C_v T + \frac{2aP}{-1+\sqrt{1+\frac{4aP}{RT}}} + q. \tag{2.26}$$

Consequently, the specific heat at constant pressure for the VO1 EOS reads,

$$C_p = C_v + \frac{4a^2P^2}{RT^2\sqrt{1+\frac{4aP}{RT}}\left(-1+\sqrt{1+\frac{4aP}{RT}}\right)^2}. \tag{2.27}$$

It is interesting to note that the specific heat at constant pressure depends on the thermodynamic state and is not a constant anymore compared to the NA EOS. A more compact form of (2.27) with the help of (2.25) is,

$$C_p = C_v + R\frac{\left(1+a\rho\right)^2}{1+2a\rho}. \tag{2.28}$$

Introducing the specific heat ratio $\gamma$ in the analogue of the Mayer's relation (2.28) leads to the following relation,

$$\gamma - 1 = \frac{R\left(1+a\rho\right)^2}{C_v\left(1+2a\rho\right)}, \tag{2.29}$$

meaning that $\gamma$ is density dependent.

Knowledge of the sound speed is needed. Equation (2.16) hands out for the VO1 EOS formulation,

$$c^2\left(P,\rho\right) = \frac{P}{\rho}\left(\frac{R}{C_v}\left(1+a\rho\right)+\frac{1+2a\rho}{1+a\rho}\right),$$
$$c^2\left(P,\rho\right) = \frac{1}{1+a\rho}\left[c_{IG}^2\left(1+2a\rho\right)+\frac{a^2\rho RP}{C_v}\right]. \tag{2.30}$$

Expression of the specific entropy is given in Appendix A. Convexity of the VO1 EOS is addressed in Appendix B.

Having in hands NA and VO1 EOS in thermodynamically compatible forms, a method to determine their parameters is examined. Such a methodology is presented hereafter.

## 3. Methodology to determine the NA and VO1 parameters

Determination of thermodynamic parameters for a particular reactive material is usually done by using either experimental data obtained from tests in closed vessels, and/or with computed data from a thermochemical code, such as CHEETAH (Fried and Souers, 1994).

CHEETAH employs the semi-empirical Becker-Kistiawosky-Wilson (BKW) EOS (Becker, 1921), (Kistiakowsky and Wilson, 1941). Related parameters determination is addressed for example in Hobbs and Baer (1993). Various thermochemical codes have been tested as alternatives, such as ICT (Volk and Bathelt, 1988) and BAGHEERA (Bac, 1984) both based on high order virial expansions. Constant volume explosion pressures computed with these codes have been compared to closed bomb vessel experiments. All computed pressures appeared significantly overestimated with errors between 10% to 20%. Such discrepancies have serious consequences on reduced EOS parameters determination, for instance on most fundamental data like condensed material energy. Error at this stage may have dramatic consequences on interior ballistics code predictions when applied to realistic configurations. Additionally, it appeared that all thermochemical codes predictions were in close agreement regarding flame temperature and specific heat ratio $\gamma = C_p/C_v$. Flame temperatures and coefficient $\gamma$ also seemed constant for high loading densities. The following method relies on these observations to determine accurate reduced EOS parameters. It is based on a method proposed by Farrar and Leeming (1983) with some modifications and additional consideration regarding non-ideal effects such as wall heat losses.

### *3.1. Determination of the NA EOS parameters*

For a constant volume explosion, EOS (2.4) and (2.5) read,

$$P_{max} = \frac{R\left(e\left(T_{flame}\right)-q\right)}{C_v\left(v-b\right)}, \tag{3.1}$$

$$e\left(T_{flame}\right) = C_v T_{flame} + q, \tag{3.2}$$

where $T_{flame}$ is the gas products flame temperature of the energetic material considered and $P_{max}$ is the corresponding peak pressure. Gas products thermodynamic parameters R, b and $C_v$ are associated to a given condensed material and must be determined. In this context, the specific volume v corresponds to the loading density of the reactant. The loading density is defined as the ratio of mass of solid reactive material by the combustion chamber's volume, $\rho_{load} = m_{reactive}/V_{chamber}$.

Assuming presence of non-ideal effects, such as wall thermal losses, or deviation from thermochemical equilibrium during the constant volume explosion, energy is not constant during the whole process,

$$e\left(T_{flame}\right) - q = e_s - \pi. \tag{3.3}$$

The left-hand side represents the gas energy while the right one corresponds to the solid energy $e_s$ reduced by the various losses $\pi$. The difference $e_s - \pi$ represents the effective energy $e_{s,eff}$ i.e. the energy received by the gas phase through solid phase combustion. Expressions (3.1) and (3.2) are rewritten as,

$$P_{max} = \frac{R T_{flame}}{\left(v-b\right)}, \tag{3.4}$$

$$e_{s,eff} = e_s - \pi = C_v T_{flame}. \tag{3.5}$$

Consequently, the effective energy $e_{s,eff}$ is a key thermodynamic parameter to determine for a given reactive material. In all thermochemical codes mentioned previously, wall heat losses $\pi$ are not considered leading to overestimated pressures compared to closed bomb vessel experiments. This feature prevents direct computation of the condensed material effective energy with a thermochemical code. Closed bomb vessel pressure measurements are preferred as they include implicitly non-ideality such as wall heat losses.

The following methodology based on closed vessel pressures is presented in Farrar and Leeming (1983). Two closed vessel tests are considered with two different masses of reactive material $m_1$ and $m_2$ and the corresponding peak pressures $P_{max,1}$ and $P_{max,2}$ are measured. The following system of two equations is then constructed,

$$\begin{cases} P_{max,1}\left(V_{chamber} - m_1 b\right) = m_1 F \\ P_{max,2}\left(V_{chamber} - m_2 b\right) = m_2 F \end{cases}. \tag{3.6}$$

Farrar and Leeming (1983) solve this system to determine the covolume b and $F = RT_{flame}$ defined as the "force" constant of the reactive material. However, in this work the method is slightly modified to determine both b and R by relying on the computed flame temperature for two main reasons. The first one is to adjust parameters at a certain flame temperature which improve EOS temperature predictions. The second reason is related to the mixture thermodynamic models developed in Section 6 which require independent knowledge of R for each gas produced by a given reactive material. This revised method is presented hereafter.

Using two points $\left(P_{max,1}, v_1\right)$ and $\left(P_{max,2}, v_2\right)$ from experimental measurements, and the flame temperature $T_{flame}$ computed with a thermochemical code, the following algebraic system is obtained,

$$\begin{cases} P_{max,1} = \dfrac{RT_{flame}}{v_1 - b} \\ P_{max,2} = \dfrac{RT_{flame}}{v_2 - b} \end{cases}. \tag{3.7}$$

Here, the flame temperature is considered constant, at least in the present high-density conditions of the two closed bomb points. The validity of this assumption will be discussed on practical examples in Section 4. Resolution of (3.7) hands out the expressions for the parameters,

$$\begin{aligned} b &= \frac{P_{max,1}v_1 - P_{max,2}v_2}{P_{max,1} - P_{max,2}}, \\ R &= T_{flame}P_{max,1}P_{max,2}\left(\frac{v_2 - v_1}{P_{max,1} - P_{max,2}}\right). \end{aligned} \tag{3.8}$$

The two points used to determine the parameters define density and pressure ranges where these parameters are valid. It is possible to use more points, if corresponding data are available, and use a fitting method to determine parameters b and R. However, such a method does not significantly improve accuracy.

The methodology to determine parameters b and R now available, determination of $C_v$ and $e_{s,eff}$ is addressed hereafter.

The specific heat ratio $\gamma$ is considered constant and computed with a thermochemical code. Indeed, after computations with the various thermochemical codes for different loading densities, coefficient $\gamma$ appeared weekly dependent on density and pressure, at least in the studied range between approximately 100 kg/m$^3$ and 400 kg/m$^3$. At this stage, the specific heat $C_v$ is also considered constant and computed with $\gamma$ and the Mayer's relation for the NA EOS (2.15),

$$C_v = \frac{R}{\gamma - 1}. \tag{3.9}$$

Constant specific heat validity will be studied deeper in Section 5 considering a temperature dependent specific heat.

The effective energy $e_{s,eff}$ is deduced with equation (3.5) and the flame temperature.

This method based on two experimental closed bomb vessel points is particularly accurate for gun's computation. But its validity is restricted to the density and pressure ranges used for parameters determination. As will be seen in Section 4, in the frame of the NA EOS, when the density comes out of this range, the computed pressure becomes inaccurate. This issue becomes pregnant for advanced ballistic systems with enhanced loadings. This is the motivation of the following investigation.

### *3.2. Determination of the VO1 EOS parameters*

For a constant volume explosion, EOS (2.19) and (2.20) read,

$$P_{max} = \rho R T_{flame} \left(1 + a\rho\right), \tag{3.10}$$

$$e_{s,eff} = e_s - \pi = C_v T_{flame}. \tag{3.11}$$

As before, the thermodynamic parameters R, $e_{s,eff}$, a and $C_v$ are associated to a given material and must be determined.

Similar approach as for the NA EOS is used to determine the specific gas constant R and the virial coefficient a. Knowing two points $(P_{max,1}, \rho_1)$ and $(P_{max,2}, \rho_2)$ from closed bomb vessel experiments and flame temperature $T_{flame}$ from thermochemical computations, the following non-linear system is expressed with equation (3.10),

$$\begin{cases} P_{max,1} = \rho_1 R T_{flame} \left(1 + a\rho_1\right) \\ P_{max,2} = \rho_2 R T_{flame} \left(1 + a\rho_2\right) \end{cases}. \tag{3.12}$$

Resolution leads to the expressions,

$$a = \frac{P_{max,2}\rho_1 - P_{max,1}\rho_2}{P_{max,1}\rho_2^2 - P_{max,2}\rho_1^2},$$
$$R = \frac{1}{T_{flame}} \left( \frac{P_{max,1}\rho_2^2 - P_{max,2}\rho_1^2}{\rho_1\rho_2^2 - \rho_1^2\rho_2} \right). \tag{3.13}$$

Assuming again constant specific heat ratio $\gamma$, parameter $C_v$ is determined using the Mayer's relation corresponding to the VO1 formulation (2.29),

$$C_v = \frac{R}{(\gamma - 1)} \frac{\left(1 + a\bar{\rho}\right)^2}{\left(1 + 2a\bar{\rho}\right)}. \tag{3.14}$$

The density used in (3.14) to compute $C_v$ is the average loading density of the two "experimental" points available, $\bar{\rho} = (\rho_1 + \rho_2)/2$. This choice is arbitrary and unimportant. One could use $\rho_1$ or $\rho_2$ or even an arbitrary density outside the adjustment range since both $\gamma$ and R are nearly constant for high loading densities. Indeed, in high density/pressure conditions, dissociation effects are limited leading to nearly constant composition of combustion gases and consequently invariant R.

Finally, the effective energy $e_{s,eff}$ is determined with equation (3.11) and the flame temperature.

Having in hands parameters for both NA and VO1 EOS, these reduced formulations are compared to the predictions given by CHEETAH using BKW EOS, considered as reference.

## 4. Application to a single reactive material

In all tests presented in this section and subsequent ones, thermodynamic parameters $R, b$ and a have been obtained using the two-points method described in Section 3. However, experimental data points $(P_{max}, \rho)_{(1,2)}$ are replaced by CHEETAH computations at two different loading densities to illustrate the method on various examples. Consequently, wall heat losses and other non-ideal effects are considered absent $(\pi = 0)$ allowing comparisons between reduced EOS (NA and VO1) and CHEETAH computations, based on BKW EOS and serving as reference. When closed bomb vessel experimental data are available, the same method is applied. With experimental data, wall heat losses are accounted for, their consideration having significant consequences on computed pressures.

The methodology described in Section 3 is applied to determine NA and VO1 parameters resulting of various reactive materials combustion. Considered materials are nitrocellulose 13% nitrogen (NC-13), RDX, nitroglycerin (NG) and HMX which are reactive materials widely used in propellant charges. Computed data with the thermochemical code CHEETAH for constant volume explosions are given in Table I.

| **Material** | $\rho_1$ $(kg/m^3)$ | $P_{max,1}$ $(MPa)$ | $\rho_2$ $(kg/m^3)$ | $P_{max,2}$ $(MPa)$ | $T_{flame}$ $(K)$ | $\gamma$ |
|---|---|---|---|---|---|---|
| NC-13 | 100 | 130.3 | 150 | 214.1 | 3275 | 1.207 |
| RDX | 100 | 163.4 | 150 | 267.6 | 4040 | 1.214 |
| NG | 100 | 131.6 | 150 | 215.1 | 3991 | 1.180 |
| HMX | 100 | 162.3 | 150 | 265.7 | 4012 | 1.211 |

**Table I.** Closed bomb vessel pressure computed with CHEETAH for various reactive material at two different loading densities. Flame temperatures and specific heat ratio are computed at an arbitrary high loading density to adjust temperature predictions, here 200 kg/m$^3$.

These data are used to determine the parameters for NA and VO1 EOS and are representative of loading densities currently used in closed bomb vessels. The parameters are thus adjusted in the density range $[100\,; 150\ kg/m^3]$.

| Material | EOS | $C_v$ (J/kg/K) | R (J/kg/K) | $e_{s,eff}$ (kJ/kg) | b \| a $(m^3/kg)$ |
|---|---|---|---|---|---|
| NC-13 | NA | 1637.1 | 338.9 | 5360.7 | 0.001484 |
| | VO1 | 1640.5 | 322.0 | 5371.9 | 0.002359 |
| RDX | NA | 1640.9 | 346.2 | 6629.3 | 0.001440 |
| | VO1 | 1644.1 | 330.2 | 6642.1 | 0.002249 |
| NG | NA | 1573.1 | 283.2 | 6277.9 | 0.001413 |
| | VO1 | 1576.0 | 270.6 | 6289.5 | 0.002185 |
| HMX | NA | 1642.0 | 346.5 | 6588.5 | 0.001435 |
| | VO1 | 1645.2 | 330.6 | 6601.1 | 0.002237 |

**Table II.** NA and VO1 EOS parameters for various reactive materials adjusted in the density range $[100\,;150\ \mathrm{kg/m^3}]$.

For each material, parameter $C_v$ is computed for the VO1 formulation with equation (3.14) at the density $\bar{\rho} = (\rho_1 + \rho_2)/2 = 125\ \mathrm{kg/m^3}$. As mentioned in the previous section, this choice is arbitrary and without consequences. Indeed, for example at $\rho_1 = 100\ \mathrm{kg/m^3}$ for gases produced by NC-13 combustion, the specific heat ratio $\gamma = 1.211$ while at $\rho_2 = 150\ \mathrm{kg/m^3}$, $\gamma = 1.209$. Such variations are insignificant when compared to the loading density variation. Similar observations are made for all considered materials.

Reference curves computed with CHEETAH and predicted results for each reduced EOS with the data of Table II are compared in Figure 1 and Figure 2. Comparisons are made for both flame temperature and maximum pressure on a wide range of loading densities. NA and VO1 having the same caloric EOS, flame temperatures for both formulations are computed as $T_{flame} = e_{s,eff}/C_v$. Pressure for NA and VO1 EOS are computed respectively with (3.4) and (3.10).

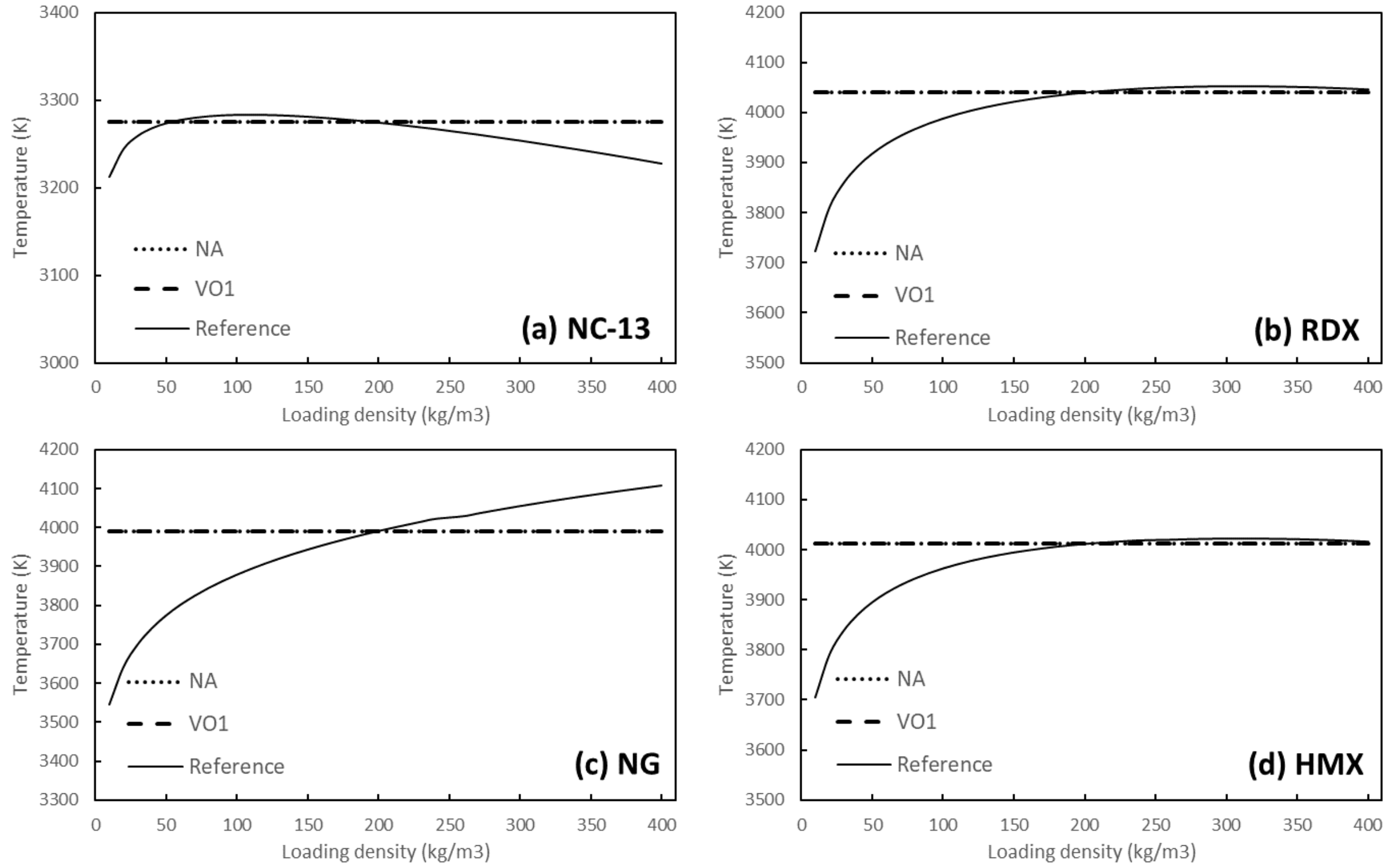

**Figure 1.** Flame's temperature of various reactive materials computed for different loading densities. Reduced NA and VO1 EOS are compared to CHEETAH computations, considered as reference. The black solid line represents the reference results. The black dotted line represents the NA EOS results. The black dashed line represents the VO1 EOS results. Both reduced EOS predictions are identical and equal to the flame temperature used in parameters determination (Table I).

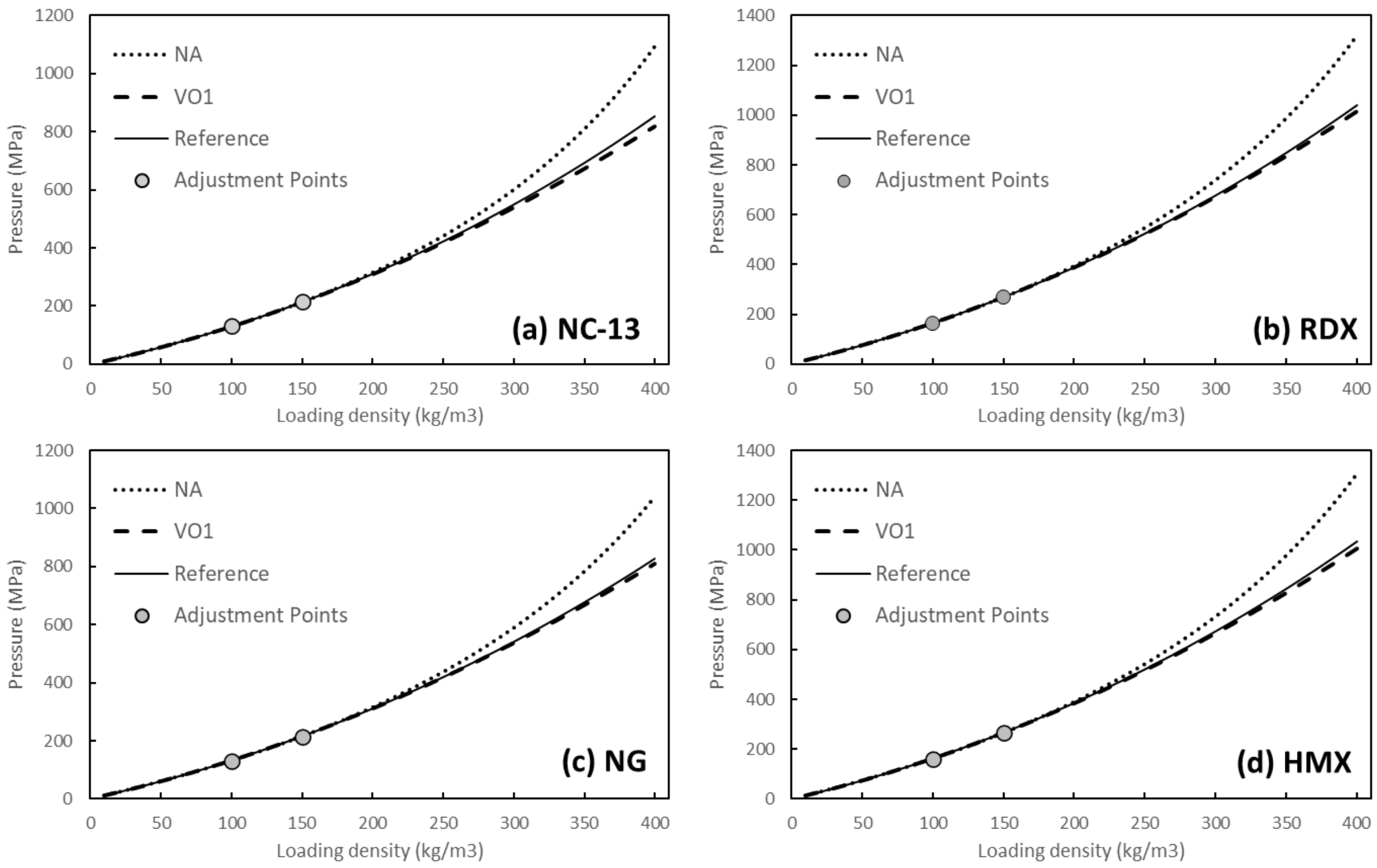

**Figure 2.** Constant volume explosion pressure of various reactive materials computed for different loading densities. The black solid line represents reference results obtained with the thermochemical code CHEETAH. The black dotted line represents the NA EOS results. The black dashed line represents the VO1 EOS results. Adjustment points indicate data used in the determination of NA parameters with (3.8) and VO1 parameters with (3.13). Extrapolation at high density clearly shows deviations resulting of NA formulation, compared to VO1. At density 400 kg/m3, NA pressure error against the reference is approximately 250 MPa (27%) for all materials while VO1 error is about 27 MPa (3%).

It appears immediately that the flame temperature is nearly constant in these computations, especially for loadings higher than 100 kg/m$^3$ as shown in Figure 1 by the reference curves. Indeed, high pressure conditions limit molecular vibration and forbid dissociation effects. Even though nitroglycerin products temperature does not seem constant, between densities 100 kg/m$^3$ and 400 kg/m$^3$, the flame temperature increases by merely 6% while the loading density is quadrupled. The flame temperature can thus be considered constant in the present applications. The same observations are made for the specific heat ratio $\gamma$ which varies by approximately 1% between 100 kg/m$^3$ and 400 kg/m$^3$ for each material.

In Figure 1, flame temperatures predicted by the reduced EOS correspond to the one use in the various parameters determination (see Table I). Consequently, predicted results are in good agreement since the temperature is only slightly varying with the loading density.

Regarding computed pressure, for the various reactive materials all three curves are in very close agreement in the range used to calibrate the EOS. For higher loading densities, deviations between reference and EOS predictions gradually appear. These deviations are much higher for NA EOS. The VO1 EOS stays quite close to the reference for every considered material. Figure 2 clearly shows what was stated in introduction and is one of the main motivations of this work. When extrapolated, the NA EOS diverges from the reference solution while the VO1 EOS remains much closer.

A last remark is that adjusting the thermodynamic parameters directly for higher densities would obviously provide better results by resolving the extrapolation problem. However, the thermochemical code CHEETAH is only used here to illustrate the two points method presented in Section 3 by providing data representative of closed bomb tests (Table I), which are often limited to “low” loadings. Once more, it is advised to use experimental closed bomb vessel data in any real application for reasons stated prior in Section 3.

As the VO1 formulation appears simple and accurate, attempt is done to extend it to temperature dependent specific heat, to try to cover a wider range of temperatures. Large temperature range appear during gun’s ignition stage. This is the aim of the following section.

## 5. Extension of the VO1 EOS with $C_v(T)$

### *5.1. Formulation with a temperature dependent specific heat*

In this section, extension of the VO1 EOS to temperature dependent specific heat is examined.

Linear dependency of the specific heat is assumed to keep simple formulations and make parameters determination easy. Thus, the specific heat reads,

$$C_v(T) = C_{v_0} + cT, \tag{5.1}$$

with $C_{v_0}$ and c constant parameters defined for a given material. The caloric EOS (2.2) now reads,

$$e(T) = \overline{C_v}(T)(T - T_0) + e_0. \tag{5.2}$$

$\overline{C_v}(T)$ is the average specific heat in the temperature interval $[T_0, T]$ and is given by,

$$\overline{C_v}(T) = \frac{1}{T - T_0} \int_{T_0}^{T} C_v(T)\,dT = C_{v_0} + \frac{c}{2}(T + T_0) \tag{5.3}$$

An energy constant q is defined to simplify EOS formulations,

$$q = e_0 - \left( C_{v_0} T_0 + \frac{c}{2} T_0^2 \right). \tag{5.4}$$

Inserting (5.3) and (5.4) in (5.2), the VO1 thermal and caloric EOS with variable specific heat read,

$$P(\rho, T) = \rho R T (1 + a\rho), \tag{5.5}$$

$$e(T) = C_{v_0} T + \frac{c}{2} T^2 + q. \tag{5.6}$$

The temperature is determined by solving the quadratic equation provided by (5.6). Only the positive root is retained,

$$T = \frac{-C_{v_0} + \sqrt{C_{v_0}^{\ 2} + 2c(e - q)}}{c}. \tag{5.7}$$

Inserting (5.7) in (5.5) results in the following function,

$$P(\rho, e) = \frac{\rho R}{c} \left[ \sqrt{C_{v_0}^{\ 2} + 2c(e - q)} - C_{v_0} \right] (1 + a\rho). \tag{5.8}$$

With this extended VO1 formulation, the internal energy (5.6) is still temperature dependent only, and the function $P = P(\rho, T)$ is identical to the one for VO1 with constant specific heat. Consequently, thermodynamic compatibility (2.23) is inherently fulfilled.

The various parameters of this thermodynamic model need to be determined. Considering again constant volume explosion tests, methods employed in Section 3 are reused to find parameters R, a and $e_{s,eff}$. However, parameters of the variable specific heat $C_{v_0}$, c and the reference energy q cannot be determined without modifications. Indeed, as observed in Section 4, the flame temperature in constant volume explosions is nearly constant for a given material at different loading densities. Therefore, to make the temperature vary, an inert chemical species is considered with the reactant. Varying the proportion of this inert will make the flame temperature vary. Consequently, mixture VO1 EOS with $C_v(T)$ is needed. This is the aim of the following derivations.

Considering gas products resulting of the reactive material combustion, with mass fraction Y, and the inert species with mass fraction of $1 - Y$, the mixture internal energy reads,

$$e_{mix}(T) = Ye(T) + (1-Y)e_{in}(T), \tag{5.9}$$

where the subscript "in" refers to the inert component.

Inserting the caloric EOS (5.6) leads to the mixture caloric EOS,

$$e_{mix}(T) = Y\left(C_{v_0}T + \frac{c}{2}T^2 + q\right) + (1-Y)\left(C_{v_0,in}T + \frac{c_{in}}{2}T^2 + q_{in}\right). \tag{5.10}$$

The mixture thermal EOS is deduced of (5.5) as,

$$P(\rho_{mix}, T) = \rho_{mix}R_{mix}T(1 + a\rho_{mix}), \tag{5.11}$$

where $\frac{1}{\rho_{mix}} = \frac{Y}{\rho} + \frac{1-Y}{\rho_{in}}$ and $R_{mix} = \hat{R}\left(\frac{Y}{\hat{W}} + \frac{1-Y}{\hat{W}_{in}}\right)$.

Having in hands the VO1 EOS with temperature dependent specific heat, parameters determination is addressed.

### *5.2. Determination of the VO1 $C_v(T)$ parameters*

Parameters for the inert component $C_{v_0,in}$, $c_{in}$ and $q_{in}$ are supposed known beforehand. More precisely, noble gases are considered as inert species (Argon or Xenon), with constant specific heat and zero heat of formation ($c_{in} = 0, q_{in} = 0$). The constant specific heat $C_{v_0,in} = C_{v_{in}}$ is considered known.

Parameters $C_{v_0}$, c and q need to be determined for gas products, resulting of the reactive material combustion. Internal energy invariance in the closed bomb vessel reads,

$$e_{mix}(T_{flame}) = Y\left(C_{v_0}T_{flame} + \frac{c}{2}T_{flame}{}^2 + q\right) + (1-Y)\left(C_{v_0,in}T_{flame} + \frac{c_{in}}{2}T_{flame}{}^2 + q_{in}\right). \tag{5.12}$$

CHEETAH is used with variable mass fraction Y making the flame temperature $T_{flame}$ vary. With these informations, a method is developed in Appendix C for the practical determination of parameters $C_{v_0}$, c and q. An important assumption is that the gas products composition is considered constant when diluted in the inert species. Indeed, the specific heat $C_v(T)$ is intended to consider temperature variations of gas but invariance in its composition

Parameters R and a are determined following the same methodology presented in Section 3. For a constant volume explosion, equations (5.5) and (5.6) become,

$$P_{max} = \rho R T_{flame}\left(1 + a\rho\right), \tag{5.13}$$

$$e_{s,eff} = e\left(T_{flame}\right) - q = C_{v_0} T_{flame} + c T_{flame}^{\ 2}. \tag{5.14}$$

Equation (5.13) is identical to equation (3.10) for VO1 EOS with constant specific heat. Consequently, knowing the numerical values at two points $\left(P_{max,1}, \rho_1\right)$ and $\left(P_{max,2}, \rho_2\right)$ from closed bomb vessel experiments data and the flame temperature with a thermochemical code, the same system (3.12) is available, leading to the same parameter expressions (3.13).

The effective energy $e_{s,eff}$ is computed with equation (5.14) and the flame temperature.

The specific gas constant $R_{mix}$ of gas products resulting of the reactive material combustion in an inert is deduced from molar masses $\widehat{W}$ and $\widehat{W}_{in}$. The reactive material molar mass is computed as $\widehat{W} = \widehat{R}/R$. For the inert gas, $\widehat{W}_{in}$ is supposed known beforehand.

The various EOS parameters are now available, and the next section presents applications of VO1 EOS with $C_v(T)$ to different examples.

### *5.3. Application to the explosion of a reactive material with an inert component*

In this section, the VO1 EOS with a temperature dependent specific heat is tested and compared to results from NA and VO1 EOS with constant specific heat and reference results computed with the thermochemical code CHEETAH. Since the goal of the extended VO1 formulation is to improve accuracy when the temperature varies, a specific example is considered. As mentioned previously, mixtures of a reactive material (NC-13) and inert component (Argon) are considered.

Application of the methods presented in Section 5.2 and Appendix C lead to the parameters summarized in Table III. Here, only NC-13 is studied as the composition of gas products resulting of its combustion is unchanged when diluted in an inert. That is not the case for other reactive materials considered in Section 4 (RDX, NG and HMX), as shown in Appendix C. Such variations in composition are mainly due to changes in chemical equilibrium induced by pressure and/or temperature variations when combustion gases are diluted with an inert gas.

| **Parameters** | **NC-13** | **Argon** |
|---|---|---|
| $C_{v_0}$ $(J/kg/K)$ | 1416.8 | 312.2 |
| $c$ $(J/kg/K^2)$ | 0.0637 | 0 |
| $\widehat{W}$ $(g/mol)$ | 25.82 | 39.95 |
| $R$ $(J/kg/K)$ | 322.0 | 208.1 |
| $e_{s,eff}$ $(kJ/kg)$ | 4980.7 | |
| $a$ $(m^3/kg)$ | 0.002359 | |

**Table III.** VO1 EOS parameters with temperature dependent specific heat. Parameters of the gas products resulting of NC-13 combustion are adjusted in the density range $[100\,;150\ kg/m^3]$ and temperature range $[1600\,;3300\ K]$. Parameters of the inert component Argon are known from the literature (Chase, 1998).

For NA and VO1 EOS with constant specific heat, thermodynamic parameters for NC-13 and RDX are given in Table II.

Figure 3 and Figure 4 compare results obtained with the various thermodynamic formulations, and the thermochemical code CHEETAH, for several mixtures of NC-13 and Argon with mass fractions of reactive material ranging from 15% to 100%. Here, the loading density corresponds to the mixture density.

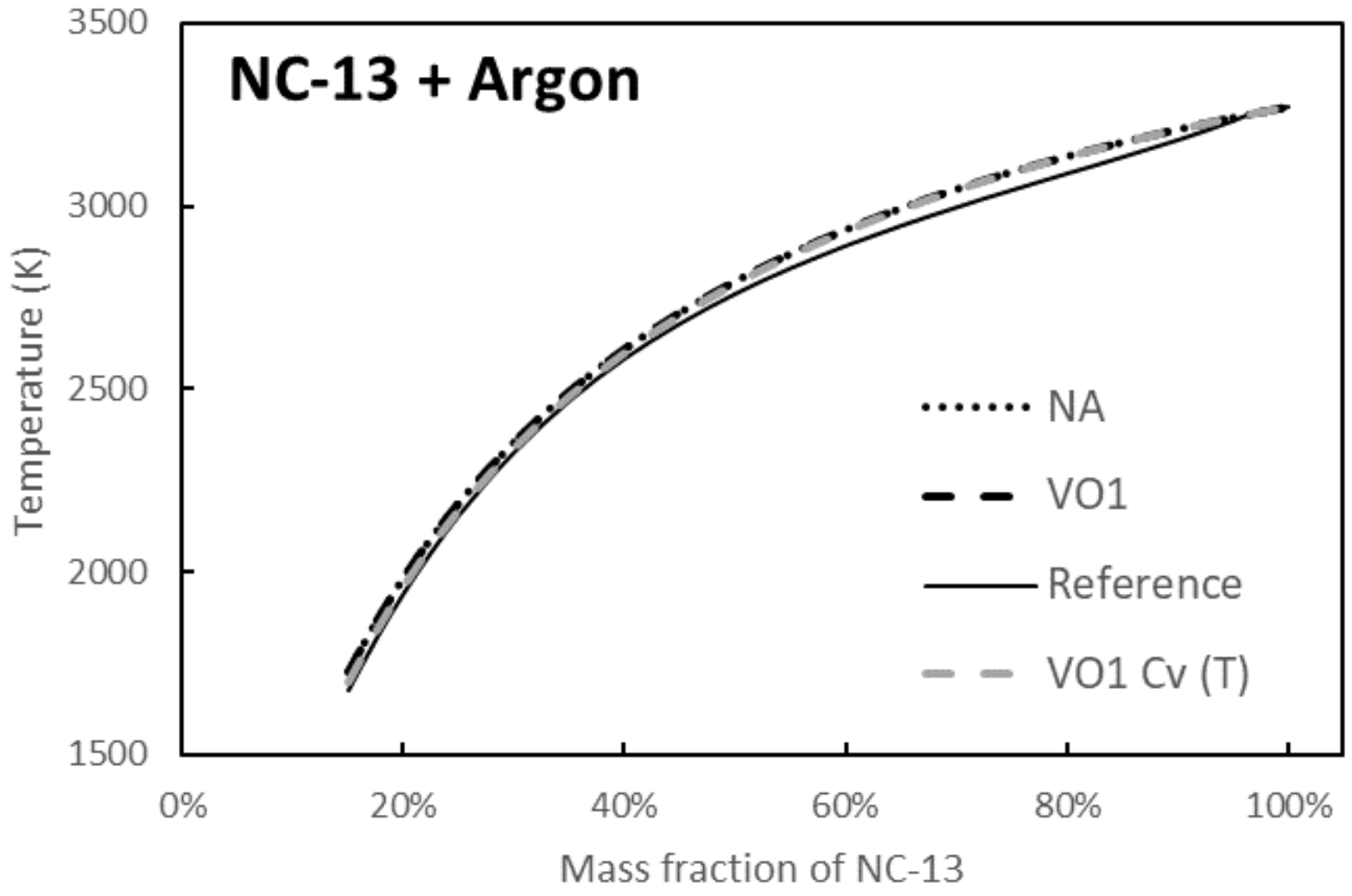

**Figure 3.** Comparison of flame temperature computation with various reduced EOS for different mixtures of NC-13 and Argon. The black solid line represents reference results obtained with the thermochemical code Cheetah. The black dotted line represents the NA EOS results. The black dashed line represents the VO1 EOS results with constant specific heat. The grey dashed line represents the VO1 EOS results with variable specific heat. All reduced EOS results appear in close agreement with the reference.

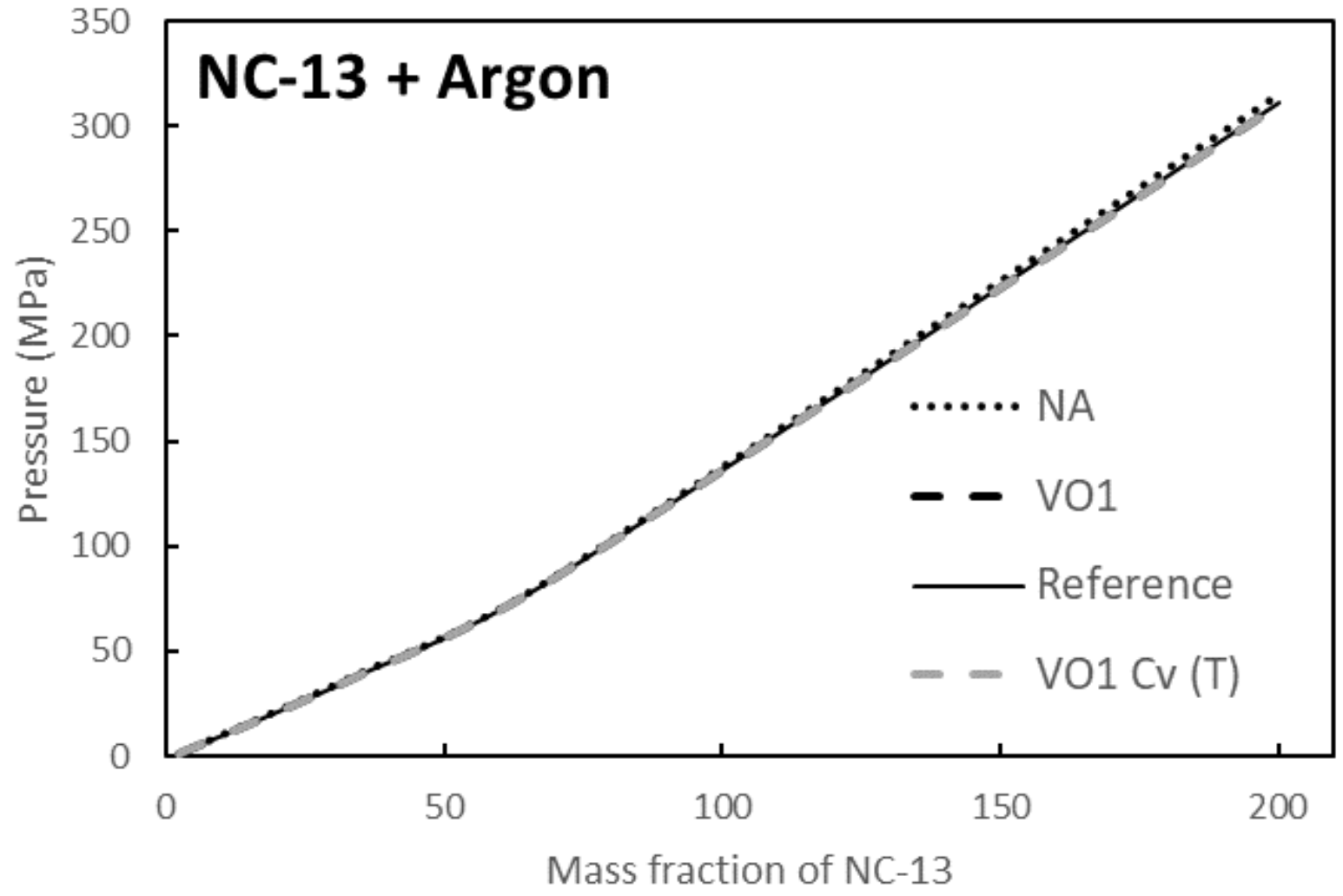

**Figure 4.** Computed pressure comparison in closed bomb vessels for different mixtures of NC-13 and Argon. The black solid line represents reference results obtained with the thermochemical code Cheetah. The black dotted line represents the NA EOS results. The black dashed line represents the VO1 EOS results with constant specific heat. The grey dashed line represents VO1 EOS results with variable specific heat. All curves are in very close agreement.

Figure 3 shows that the extended VO1 formulation with variable specific heat does not significantly improve flame temperature computation. Constant specific heat formulation seems appropriate even when the temperature varies in the range $[1600\,;3300\,\mathrm{K}]$. Figure 4 shows that all thermodynamic models are quite accurate regarding pressure computations. It is worth to note that the present density range is "low" in the sense that the loading density used in computations is not extrapolated at high levels, as was done in Figure 2. For such extrapolation at higher loading densities, VO1 EOS with constant or variable specific heat are still in perfect agreement.

To summarize, considering temperature dependent specific heat in the formulation does not improve the predictions, at least in the conditions of the present study and underlying applications. Temperature variations seems more likely dependent on chemical equilibrium changes. Such phenomena are reliant on thermochemical considerations and falls outside the scope of the present study which focuses on developing reduced thermodynamic closure. Therefore, extension of the VO1 formulation with a temperature dependent specific heat is not investigated further in the rest of this work.

The last part of this study focuses on the main goal of the present work which is the formulation of an equation of state for gas mixtures resulting of various energetic materials combustion.

## 6. Formulation of NA and VO1 EOS for multiple reactive materials

In this section, extension of both NA and VO1 EOS for mixtures of multiple reactive materials is addressed.

The specific heat of the various gas products is considered constant, and all parameters are supposed individually known for gases resulting from the combustion of each individual reactive material.

In the present mixture model, the various reactive materials for a given mixture are assumed having all the same oxygen balance sign. This condition guarantees that no significant post-combustion happens with the various gas mixtures. The present thermodynamic mixture model is not intended to account for post-combustion effects.

Each material is present with mass fraction $Y_k$ such that,

$$\sum_{k=1}^{N} Y_k = 1 \tag{6.1}$$

with subscript "k" denoting a particular reactant and N being the number of energetic materials.

Mixture equations of state that follow assume ideally mixed gases, evolving in both temperature and pressure equilibrium:

$$\begin{aligned} T_1 &= \dots = T_N = T, \\ P_1 &= \dots = P_N = P. \end{aligned} \tag{6.2}$$

Chiapolino et al. (2017) have shown that in the context of ideal gas mixtures, the well-known Dalton's law and partial pressures were consequences of these equilibrium conditions, along with mixture energy and mixture specific volume definitions used hereafter.

### *6.1. Mixture Noble-Abel (MNA) EOS*

The internal energy of the mixture reads,

$$e_{mix}(T) = \sum_{k=1}^{N} Y_k e_k(T). \tag{6.3}$$

Caloric EOS (2.5) is inserted,

$$e_{mix}(T) = \sum_{k=1}^{N} Y_k \left( C_{v_k} T + q_k \right). \tag{6.4}$$

It can be reformulated as,

$$e_{mix}(T) = C_{v,mix}T + q_{mix}, \tag{6.5}$$

with $C_{v,mix} = \sum_{k=1}^{N} Y_k C_{v_k}, \qquad q_{mix} = \sum_{k=1}^{N} Y_k q_k.$

To determine the thermal MNA EOS, the thermal EOS (2.4) for a given material k is considered,

$$P = \frac{R_k (e_k - q_k)}{C_{v_k} (v_k - b_k)}. \tag{6.6}$$

Inverting this equation, the specific volume reads,

$$v_k = \frac{R_k (e_k - q_k)}{P C_{v_k}} + b_k. \tag{6.7}$$

Inserting the caloric EOS (2.5) it becomes,

$$v_k = \frac{R_k T}{P} + b_k. \tag{6.8}$$

The mixture specific volume definition is now considered,

$$v_{mix} = \sum_{k=1}^{N} Y_k v_k (P, T). \tag{6.9}$$

Inserting (6.8),

$$v_{mix} = \sum_{k=1}^{N} Y_k \left( \frac{R_k T}{P} + b_k \right). \tag{6.10}$$

After some manipulations the pressure is obtained as,

$$P = \frac{T \sum_{k=1}^{N} Y_k R_k}{v_{mix} - \sum_{k=1}^{N} Y_k b_k}. \tag{6.11}$$

Replacing the temperature with the help of (6.4) finally yields the thermal MNA EOS,

$$P(v_{mix}, e_{mix}) = \frac{\sum_{k=1}^{N} Y_k R_k \left( e_{mix} - \sum_{k=1}^{N} Y_k q_k \right)}{\sum_{k=1}^{N} Y_k C_{v_k} \left( v_{mix} - \sum_{k=1}^{N} Y_k b_k \right)}. \tag{6.12}$$

In compact form, the MNA EOS (6.12) reads,

$$P(v_{mix}, e_{mix}) = \frac{R_{mix} (e_{mix} - q_{mix})}{C_{v,mix} (v_{mix} - b_{mix})}, \tag{6.13}$$

with $R_{mix} = \sum_{k=1}^{N} Y_k R_k = \hat{R} \sum_{k=1}^{N} \frac{Y_k}{\hat{W}_k} = \frac{\hat{R}}{\hat{W}_{mix}},$

$q_{mix} = \sum_{k=1}^{N} Y_k q_k, \qquad C_{v,mix} = \sum_{k=1}^{N} Y_k C_{v_k}, \qquad b_{mix} = \sum_{k=1}^{N} Y_k b_k.$

Knowledge of the sound speed is mandatory for compressible flow models. Considering frozen mixture composition, the sound speed for a mixture of N gases following the MNA equation of state is given by (see Appendix D),

$$c_{mix}^2 = \left(\frac{\partial P}{\partial v_{mix}}\right)_{s,Y_1,...,Y_N} = \frac{v_{mix}^2\left[P + \left(\frac{\partial e_{mix}}{\partial v_{mix}}\right)_{P,Y_1,...,Y_N}\right]}{\left(\frac{\partial e_{mix}}{\partial P}\right)_{v_{mix},Y_1,...,Y_N}}. \quad (6.14)$$

The following internal energy of the mixture is deduced from the thermal MNA EOS (6.12),

$$e_{mix}\left(P, v_{mix}, Y_1,..., Y_N\right) = P\frac{\sum_{k=1}^{N} Y_k C_{v_k}}{\sum_{k=1}^{N} Y_k R_k}\left(v_{mix} - \sum_{k=1}^{N} Y_k b_k\right) + \sum_{k=1}^{N} Y_k q_k. \quad (6.15)$$

Which readily leads to the sound speed for MNA EOS,

$$c_{mix}^2 = \frac{Pv_{mix}^2\left(\frac{\sum_{k=1}^{N} Y_k R_k}{\sum_{k=1}^{N} Y_k C_{v_k}} + 1\right)}{1 - \sum_{k=1}^{N} Y_k b_k} = \frac{Pv_{mix}\left(\frac{R_{mix}}{C_{v,mix}} + 1\right)}{1 - \frac{b_{mix}}{v_{mix}}}. \quad (6.16)$$

### *6.2. Mixture first-order virial (MVO1) EOS*

NA and VO1 obey to the same caloric EOS (2.5) and (2.20). Therefore, following (6.3) through (6.5), MVO1 caloric EOS is obtained as:

$$e_{mix}\left(T\right) = C_{v,mix}T + q_{mix}. \quad (6.17)$$

The starting point to derive MVO1 thermal EOS is equation (2.22) for the gas phase associated to reactive material k,

$$P = \frac{\rho_k R_k}{C_{v_k}}\left(e_k - q_k\right)\left(1 + a_k\rho_k\right). \quad (6.18)$$

As it is a quadratic function with respect to density $\rho_k$, the latter expresses:

$$\rho_k = \frac{-1 + \sqrt{1 + \frac{4a_k C_{v_k} P}{R_k\left(e_k - q_k\right)}}}{2a_k}. \quad (6.19)$$

Using once more the mixture specific volume definition, under the form,

$$\frac{1}{\rho_{mix}} = \sum_{k=1}^{N} \frac{Y_k}{\rho_k\left(P, T\right)}, \quad (6.20)$$

an implicit function appears for the pressure computation,

$$\frac{1}{\rho_{mix}} = \sum_{k=1}^{N} \frac{2a_k Y_k}{-1 + \sqrt{1 + \frac{4a_k C_{v_k} P}{R_k \left(e_k - q_k\right)}}}. \tag{6.21}$$

This equation represents the MVO1 thermal EOS. Equation (6.21) needs appropriate iterative method for the pressure computation.

It is worth to mention that the mixture EOS (6.21) may result in costly computations when used in hyperbolic multiphase flow codes. Indeed, in these codes an iterative method is already present for pressures relaxation among the phases (Lallemand and Saurel, 2000), (Lallemand et al., 2005). Using (6.21) means two imbricated Newton loops, resulting in computational cost increase compared to the MNA alternative, where a single Newton iteration is needed in the relaxation process.

Again, the sound speed for the MVO1 EOS is needed. Since the thermal MVO1 EOS is an implicit function, the same methodology used to determine the sound speed of MNA EOS cannot be carried out here. The following definition of the sound speed, given in Appendix D, is considered,

$$c_{mix}^2 = \frac{C_{p,mix}}{C_{v,mix}} \left( \frac{\partial P}{\partial \rho_{mix}} \right)_{T, Y_1, \ldots, Y_N}, \tag{6.22}$$

with $C_{p,mix} = \sum_{k=1}^{N} Y_k C_{p_k}$.

Rearranging (6.22),

$$c_{mix}^2 = -\frac{C_{p,mix} v_{mix}^2}{C_{v,mix} \left( \frac{\partial v_{mix}}{\partial P} \right)_{T, Y_1, \ldots, Y_N}}. \tag{6.23}$$

The partial derivative present in (6.23) is obtained from (6.21) and (2.20) as:

$$\left( \frac{\partial v_{mix}}{\partial P} \right)_{T, Y_1, \ldots, Y_N} = \sum_{k=1}^{N} \frac{-4a_k^2 Y_k}{R_k T^2 \sqrt{1 + \frac{4a_k P}{R_k T}} \left( -1 + \sqrt{1 + \frac{4a_k P}{R_k T}} \right)^2}. \tag{6.24}$$

Inserting (6.24) in (6.23) yields the sound speed for the MVO1 EOS formulation,

$$c_{mix}^2 = \frac{C_{p,mix} v_{mix}^2}{C_{v,mix} \sum_{k=1}^{N} \frac{4a_k^2 Y_k}{R_k T^2 \sqrt{1 + \frac{4a_k P}{R_k T}} \left( -1 + \sqrt{1 + \frac{4a_k P}{R_k T}} \right)^2}}. \tag{6.25}$$

A more compact form is available by introducing individual gas product densities,

$$c_{mix}^2 = \frac{C_{p,mix} P}{C_{v,mix} \rho_{mix}^2 \sum_{k=1}^{N} \frac{Y_k \left(1 + a\rho_k\right)}{\rho_k \left(1 + 2a_k \rho_k\right)}}. \tag{6.26}$$

MNA and MVO1 EOS are now available. Since the various parameters are already known thanks to the methods presented in Section 3, the next section addresses verifications and validations on relevant examples.

## 7. Application for multiple reactive materials

Mixture EOS derived previously are used for the computation of constant volume explosion thermodynamic state with mixtures resulting of NC-13 and RDX combustion and mixtures resulting of NC-13 and HMX combustion. NC-13, RDX and HMX are all negative oxygen balance energetic materials. NG is absent of this section since it has a positive oxygen balance. Therefore, post combustion effects between the various gases may occur, rendering the present mixture EOS irrelevant.

Thermodynamic parameters for all reactive materials have already been determined in the density range $[100 \,; 150 \text{ kg/m}^3]$ and are reported in Table II.

Various binary mixtures with variable reactants mass fractions are computed. The mass fractions of RDX and HMX are ranging from 0% to 50%. The pressure is computed at various loading densities to estimate the impact of extrapolations outside the adjustment density range. The flame temperature is nearly independent of loading density, therefore it is only examined at the loading density 200 kg/m$^3$ in Figure 5 and Figure 7. Reference curves extracted from CHEETAH and predicted results of each mixture EOS are compared for mixtures of NC-13 and RDX in Figure 5 and Figure 6 and for mixtures of NC-13 and HMX in Figure 7 and Figure 8. Flame temperatures for both formulations are computed as $T_{flame} = e_{s,eff,mix}/C_{v,mix}$ with $e_{s,eff,mix} = \sum_{k=1}^{N} Y_k e_{s,eff,k}$. Pressures for MNA EOS are obtained with (6.13) and for MVO1 EOS by solving equation (6.21).

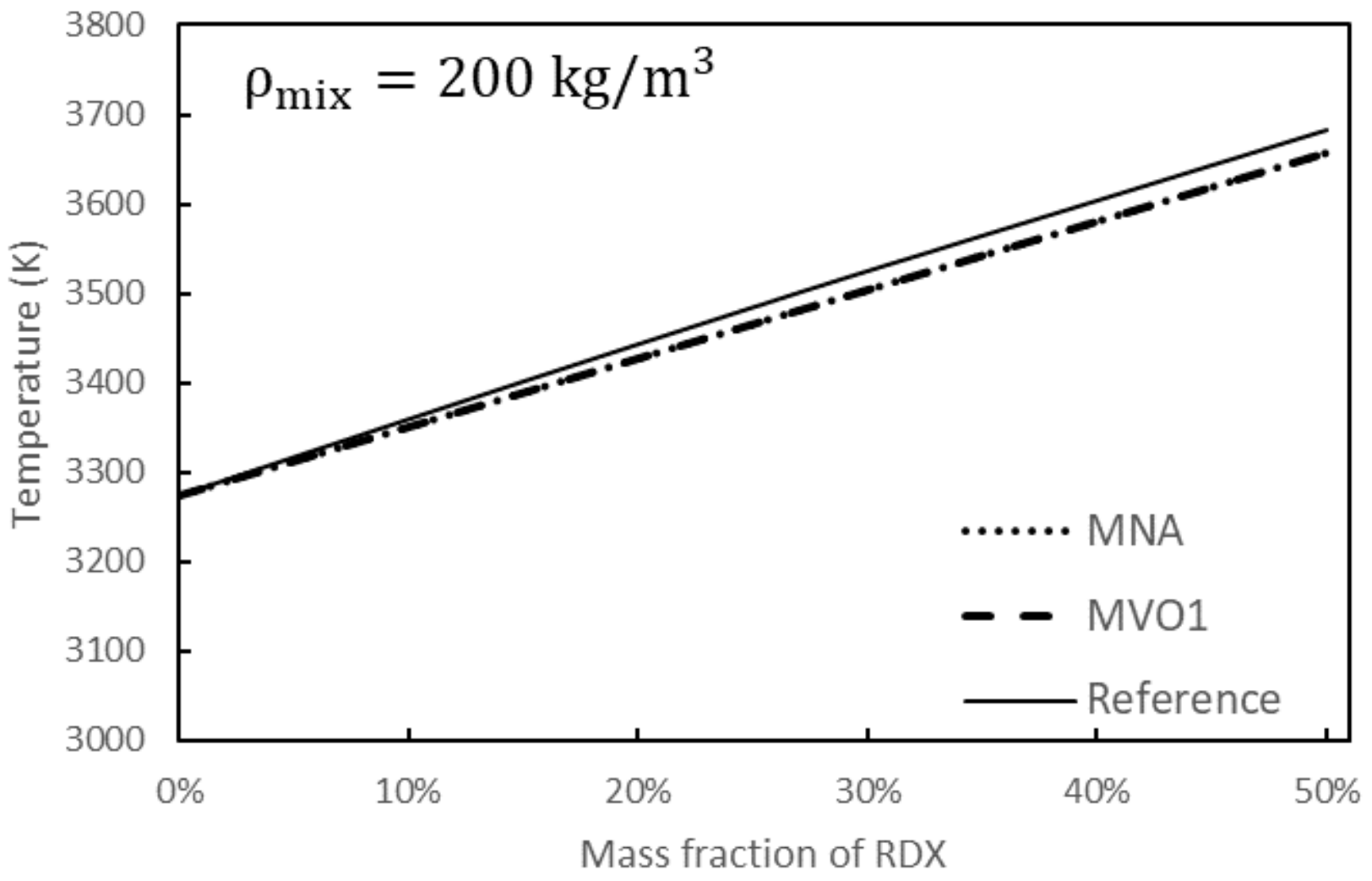

**Figure 5.** Flame temperature for various mixtures made of NC-13 and RDX at loading density 200 kg/m$^3$. The black solid line represents reference results obtained with the thermochemical code CHEETAH. The black dotted line represents the MNA EOS results. The black dashed line represents the MVO1 EOS results. Predictions for both reduced mixture EOS are identical and very close to the reference.

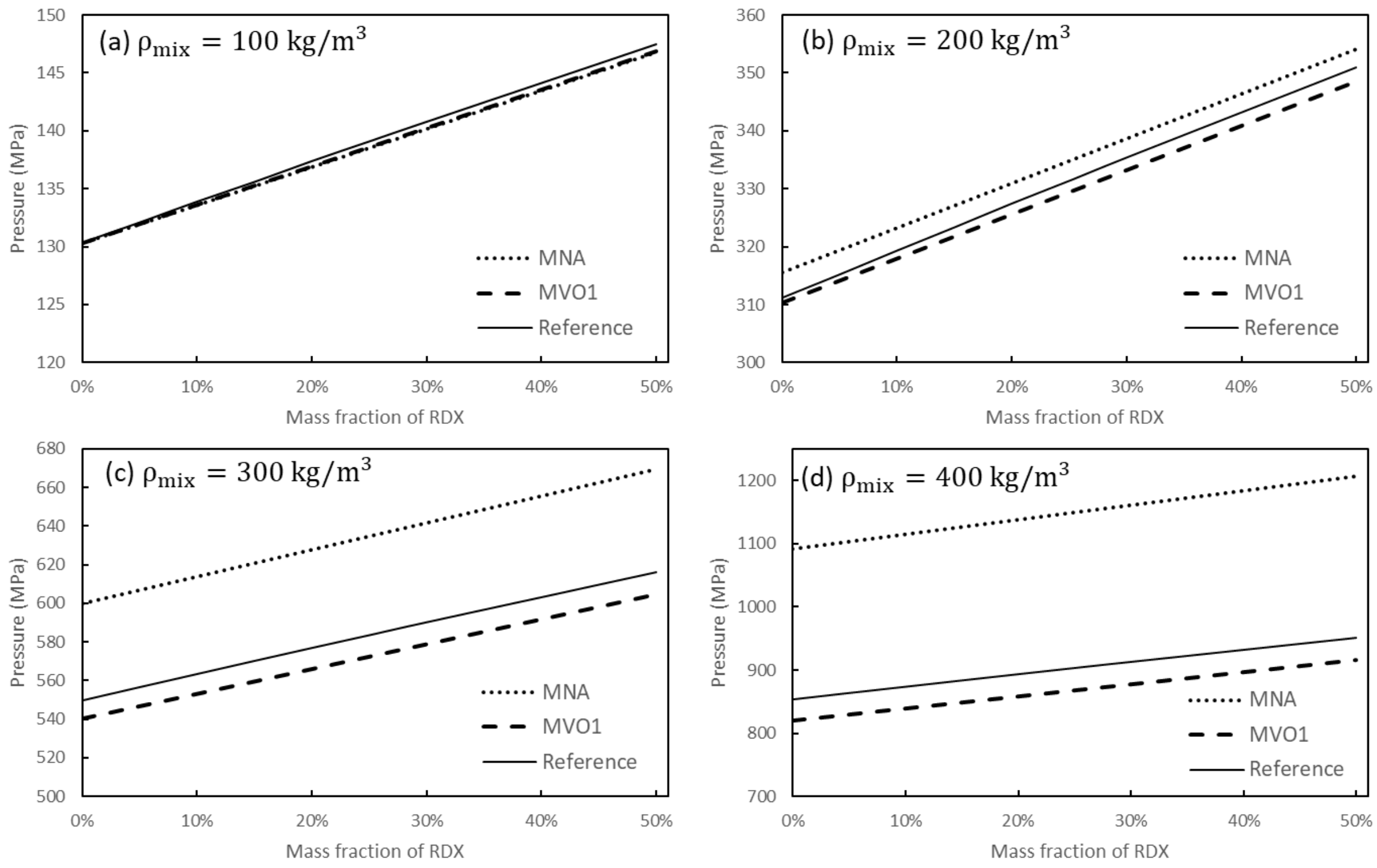

**Figure 6.** Closed bomb vessel pressure for various mixtures of NC-13 and RDX. The black solid line represents reference results obtained with the thermochemical code CHEETAH. The black dotted line represents the MNA EOS results. The black dashed line represents the MVO1 EOS results. Pressures are computed for different loading densities $\{100 ; 200 ; 300 ; 400\}$ kg/m$^3$. As before, better agreement is obtained with MVO1 than MNA, especially when extrapolating at higher densities.

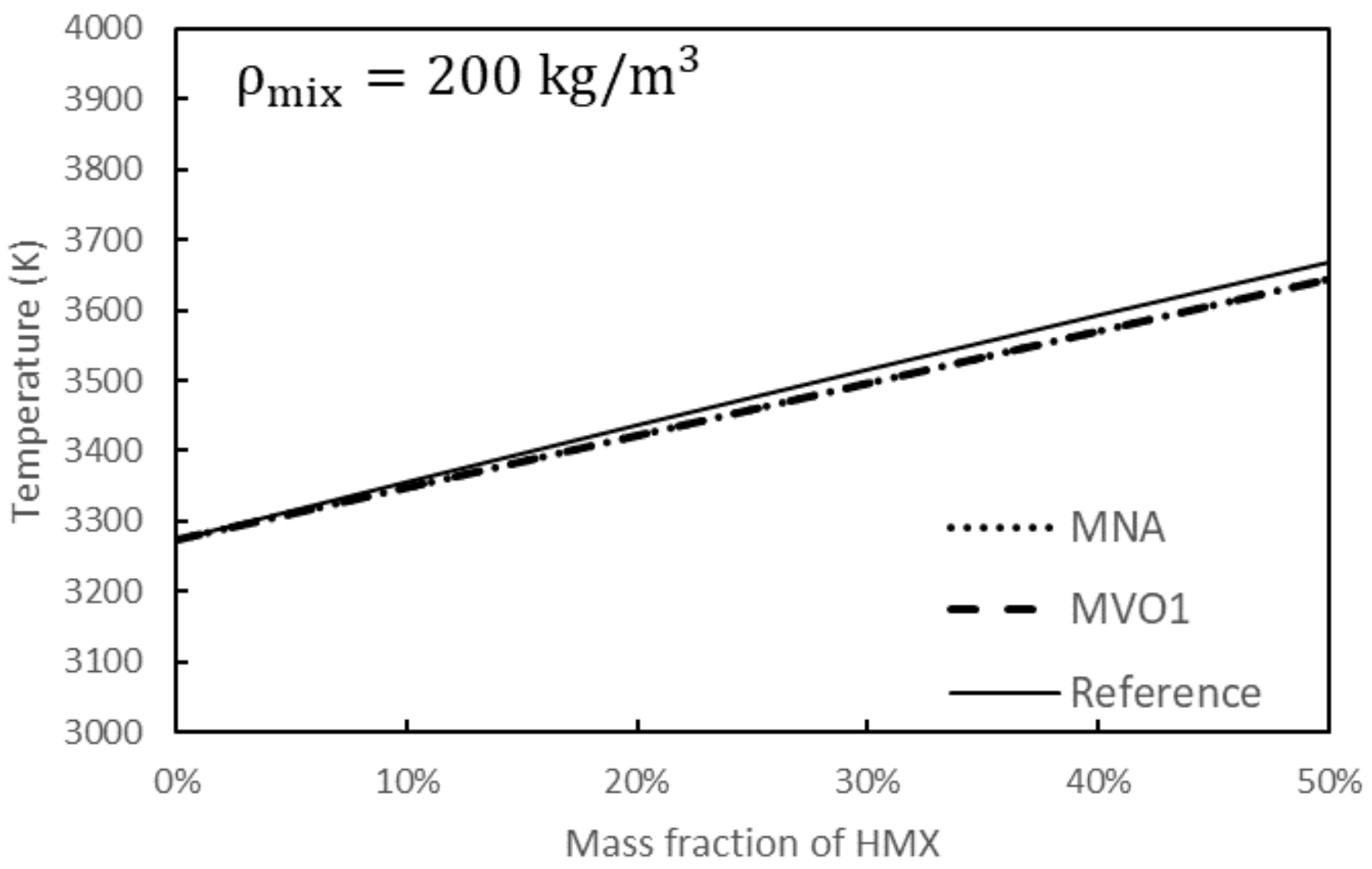

**Figure 7.** Flame temperature for various mixtures made of NC-13 and HMX at loading density 200 kg/m$^3$. The black solid line represents reference results obtained with the thermochemical code CHEETAH. The black dotted line represents the MNA EOS results. The black dashed line represents the MVO1 EOS results. Predictions for both reduced mixture EOS are identical and very close to the reference.

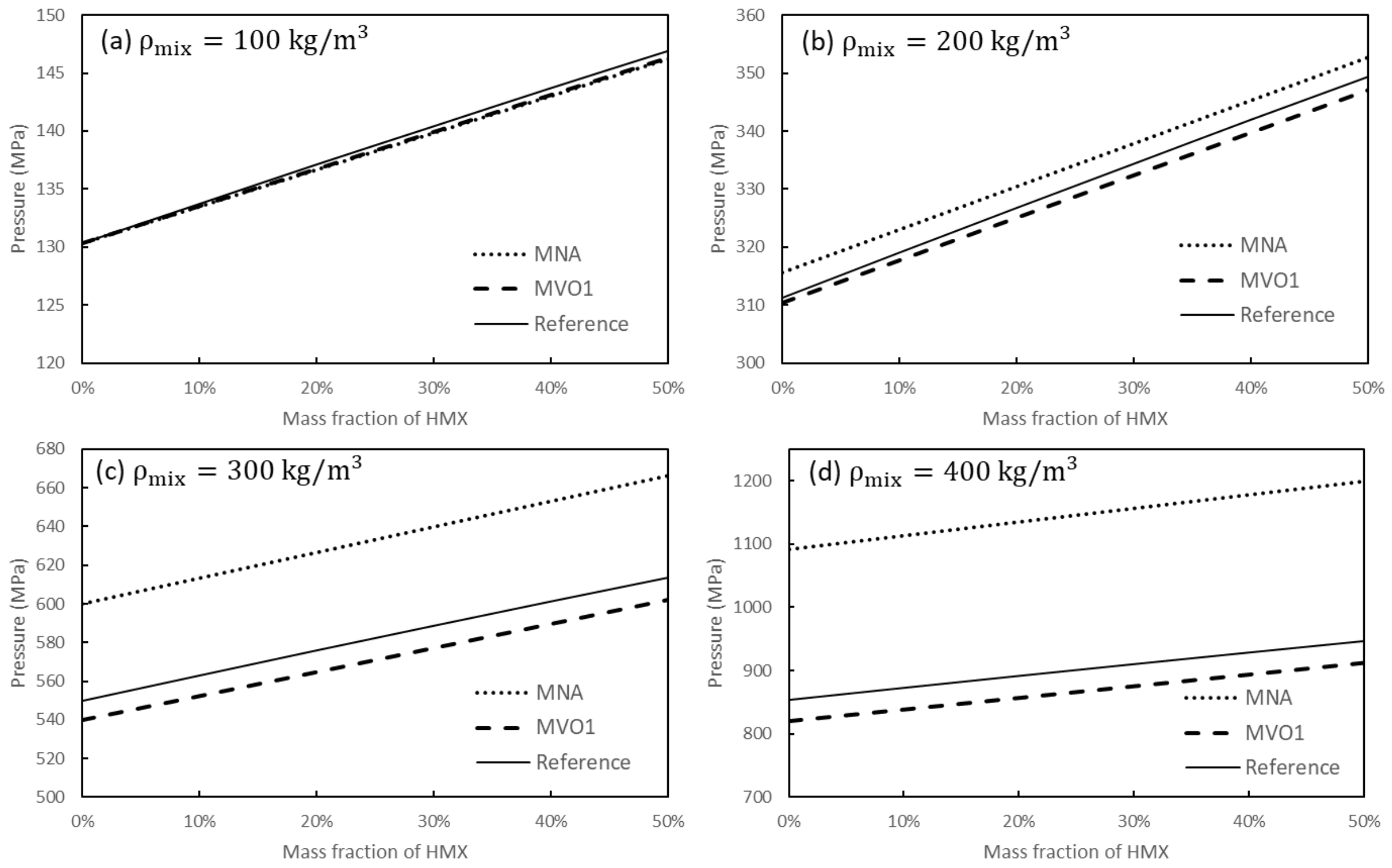

**Figure 8.** Closed bomb vessel pressure for various mixtures of NC-13 and HMX. The black solid line represents reference results obtained with the thermochemical code CHEETAH. The black dotted line represents the MNA EOS results. The black dashed line represents the MVO1 EOS results. Pressures are computed for different loading densities $\{100\,;200\,;300\,;400\}\ \mathrm{kg/m^3}$. As before, better agreement is obtained with MVO1 than MNA, especially when extrapolating at higher densities.

Regarding flame temperatures, both reduced EOS seem to predict very accurate results.

Computed pressure follows the same observations as when a single reactant was considered. Results at loading density 100 kg/m$^3$, i.e. inside the density adjustment range, are close to the reference, and gradually deviates when extrapolating at higher densities, with excessive deviations for the MNA EOS. In both examples considered, errors at loading density 400 kg/m$^3$ for MNA EOS predictions reach 250 MPa (27%) while MVO1 EOS results are much closer to the reference (35 MPa or 4%).

Overall, results given by the MVO1 EOS are more accurate than with the MNA EOS. Nonetheless, it is important to keep in mind that the MNA formulation is explicit while the MVO1 requires iterative solving for the pressure.

## 8. Conclusion

First-order virial expansion EOS has been examined as an alternative to Noble-Abel EOS to model the thermodynamic of combustion gases at high pressure and high density. The VO1 EOS appeared significantly more accurate than NA EOS for densities extrapolated outside the range used to adjust EOS parameters.

VO1 EOS being attractive, temperature dependent specific heat effects have been examined to enlarge its domain of validity. For high pressure and high-density conditions, such sophistication did not improve accuracy, when compared to predictions obtained with constant specific heat.

The NA and VO1 EOS have then been extended to deal with mixtures of materials, such situation being very common in industrial systems. Mixture NA and VO1 EOS have been developed successfully for mixtures made of reactive materials having the same oxygen balance sign. MVO1 appeared more accurate than MNA when compared to reference data. However, MNA is explicit, and consequently computationally fast. MVO1 is iterative, thus less efficient for multiphase flow computations, but more accurate, especially when the thermodynamic conditions are far from parameter's adjustment range.

## Appendix A. Entropy for the NA and VO1 EOS

Specific entropy expressions for NA and VO1 EOS are derived in this Appendix.

*A.1. Expression of the entropy for NA EOS*

The entropy formulation must satisfy the following Maxwell's relation,

$$\left(\frac{\partial s}{\partial P}\right)_T = -\left(\frac{\partial v}{\partial T}\right)_P. \tag{A.1}$$

With the help of (2.12) equation (A.1) becomes,

$$\left(\frac{\partial s}{\partial P}\right)_T = -\frac{R}{P}. \tag{A.2}$$

Integration of (A.2) is immediate,

$$s(P,T) = -R\ln(P) + K(T), \tag{A.3}$$

where $K(T)$ is a temperature dependent function. Equation (A.3) is reformulated with the specific volume,

$$s(v,T) = -R\ln\left(\frac{RT}{v-b}\right) + K(T). \tag{A.4}$$

This equation admits the partial derivative,

$$\left(\frac{\partial s}{\partial T}\right)_v = -\frac{R}{T} + \frac{dK(T)}{dT}. \tag{A.5}$$

Specific heat definition at constant volume is expressed under the following form,

$$\left(\frac{\partial s}{\partial T}\right)_v = \frac{1}{T}\left(\frac{\partial e}{\partial T}\right)_v = \frac{C_v}{T}. \tag{A.6}$$

From equations (A.5) and (A.6) the next expression arises,

$$dK(T) = (C_v + R)\frac{dT}{T}, \tag{A.7}$$

which is directly integrated to hands out the temperature dependent function,

$$K(T) = (C_v + R)\ln(T) + q', \tag{A.8}$$

with $q'$ a reference entropy. Finally inserting (A.8) in equation (A.3) yields the specific entropy expression for the NA EOS,

$$s(P,T) = -R\ln(P) + (C_v + R)\ln(T) + q'. \tag{A.9}$$

Coefficient $q'$ is defined so that the specific entropy is equal to $s_0$ at a given reference state $(P_0, T_0)$,

$$q' = s_0 + R\ln(P_0) - (C_v + R)\ln(T_0). \tag{A.10}$$

*A.2. Expression of the entropy for VO1 EOS*

The same methodology is applied to determine the specific entropy expression for the VO1 EOS. The entropy formulation must satisfy the Maxwell's relation (A.1) reformulated here as,

$$\left(\frac{\partial s}{\partial P}\right)_T = \frac{1}{\rho^2}\left(\frac{\partial \rho}{\partial T}\right)_P . \tag{A.11}$$

With the help of (2.25) equation (A.11) now reads,

$$\left(\frac{\partial s}{\partial P}\right)_T = -\frac{4a^2P}{RT^2\sqrt{1+\frac{4aP}{RT}}\left(-1+\sqrt{1+\frac{4aP}{RT}}\right)^2} . \tag{A.12}$$

Integration of (A.12) yields,

$$s(P,T) = -\frac{R}{2}\left[\sqrt{1+\frac{4aP}{RT}} + \ln\left(\left(-1+\sqrt{1+\frac{4aP}{RT}}\right)^2\right)\right] + K(T), \tag{A.13}$$

where K(T) is a temperature dependent function. Equation (A.13) is reformulated using (2.25) to introduce the density,

$$s(\rho,T) = -\frac{R}{2}\left[(1+2a\rho) + \ln\left(4a^2\rho^2\right)\right] + K(T), \tag{A.14}$$

This last equation admits the partial derivative,

$$\left(\frac{\partial s}{\partial T}\right)_\rho = \frac{dK(T)}{dT}. \tag{A.15}$$

With the specific heat definition (A.6), equation (A.15) becomes,

$$dK(T) = C_v \frac{dT}{T}, \tag{A.16}$$

which is directly integrated to hands out the temperature dependent function,

$$K(T) = C_v \ln(T) + q'. \tag{A.17}$$

Finally inserting (A.17) in equation (A.13) yields the specific entropy expression for the VO1 EOS,

$$s(P,T) = -\frac{R}{2}\left[\sqrt{1+\frac{4aP}{RT}} + 2\ln\left(-1+\sqrt{1+\frac{4aP}{RT}}\right)\right] + C_v \ln(T) + q'. \tag{A.18}$$

Coefficient $q'$ is defined so that the specific entropy is equal to $s_0$ at a given reference state $(P_0, T_0)$,

$$q' = s_0 + \frac{R}{2}\left[\sqrt{1+\frac{4aP_0}{RT_0}} + 2\ln\left(-1+\sqrt{1+\frac{4aP_0}{RT_0}}\right)\right] - C_v \ln(T_0). \tag{A.19}$$

## Appendix B. Convexity of the NA and VO1 EOS

Convexity of NA and VO1 EOS is examined. Convexity is of paramount importance for hyperbolic flow models, since closely related to wave's propagation. Convexity of an equation of state requires fulfillment of four different criteria (Godunov et al., 1976), (Menikoff and Plohr, 1989),

$$
\begin{cases}
\text{(a)} & \left(\frac{\partial^2 e}{\partial v^2}\right)_s > 0, \\
\text{(b)} & \left(\frac{\partial^2 e}{\partial s^2}\right)_v > 0, \\
\text{(c)} & \left(\frac{\partial}{\partial s}\left(\frac{\partial e}{\partial v}\right)_s\right)_v < 0, \\
\text{(d)} & \left(\frac{\partial^2 e}{\partial s^2}\right)_v \left(\frac{\partial^2 e}{\partial v^2}\right)_s - \left[\left(\frac{\partial}{\partial s}\left(\frac{\partial e}{\partial v}\right)_s\right)\right]^2 > 0,
\end{cases}
\tag{B.1}
$$

where s represents the specific entropy. However, formulation of the entropy may not always be practical to manipulate (see Appendix A). To avoid such complexity, criteria (B.1) are reformulated without entropy. Calculations are based on the two following relations.

- The Gibbs identity,

$$
de = Tds - Pdv. \tag{B.2}
$$

- Maxwell's relations (Callen and Kestin, 1960),

$$
\begin{cases}
\left(\frac{\partial T}{\partial v}\right)_s = -\left(\frac{\partial P}{\partial s}\right)_v, \\
\left(\frac{\partial T}{\partial P}\right)_s = \left(\frac{\partial v}{\partial s}\right)_P, \\
\left(\frac{\partial P}{\partial T}\right)_v = \left(\frac{\partial s}{\partial v}\right)_T, \\
\left(\frac{\partial s}{\partial P}\right)_T = -\left(\frac{\partial v}{\partial T}\right)_P.
\end{cases}
\tag{B.3}
$$

After various manipulations using (B.2) and (B.3), convexity conditions (B.1) are rearranged as:

$$
\begin{cases}
\text{(a)} & \frac{c^2}{v^2} > 0, \\
\text{(b)} & \left(\frac{\partial T}{\partial P}\right)_v \frac{\left(\frac{\partial e}{\partial v}\right)_T + P}{\left(\frac{\partial e}{\partial T}\right)_v} > 0, \\
\text{(c)} & -\frac{\left(\frac{\partial e}{\partial v}\right)_T + P}{\left(\frac{\partial e}{\partial T}\right)_v} < 0, \\
\text{(d)} & \frac{\left(\frac{\partial e}{\partial v}\right)_T + P}{\left[\left(\frac{\partial e}{\partial T}\right)_v\right]^2}\left[\left(\frac{\partial e}{\partial P}\right)_v \frac{c^2}{v^2} - \left(\frac{\partial e}{\partial v}\right)_T - P\right] > 0.
\end{cases}
\tag{B.4}
$$

These expressions are now expressed in terms of density instead of specific volume,

$$\begin{cases} \text{(a)} & \rho^2 c^2 > 0, \\ \text{(b)} & \left(\dfrac{\partial T}{\partial P}\right)_\rho \dfrac{P - \rho^2 \left(\dfrac{\partial e}{\partial \rho}\right)_T}{\left(\dfrac{\partial e}{\partial T}\right)_\rho} > 0, \\ \text{(c)} & \dfrac{\rho^2 \left(\dfrac{\partial e}{\partial \rho}\right)_T - P}{\left(\dfrac{\partial e}{\partial T}\right)_\rho} < 0, \\ \text{(d)} & \dfrac{P - \rho^2 \left(\dfrac{\partial e}{\partial \rho}\right)_T}{\left[\left(\dfrac{\partial e}{\partial T}\right)_\rho\right]^2} \left[\left(\dfrac{\partial e}{\partial P}\right)_\rho \rho^2 c^2 + \rho^2 \left(\dfrac{\partial e}{\partial \rho}\right)_T - P\right] > 0. \end{cases} \tag{B.5}$$

For the NA EOS, convexity criteria (B.4) read,

$$\begin{cases} \text{(a)} & \dfrac{c^2}{v^2} > 0, \\ \text{(b)} & \dfrac{P(v-b)}{RC_v} > 0, \\ \text{(c)} & -\dfrac{P}{C_v} < 0, \\ \text{(d)} & \dfrac{P^2}{RC_v} > 0. \end{cases} \tag{B.6}$$

As $P, T, v, R$, and $C_v$ are positive, the NA EOS is convex under the following condition:

$$b < v. \tag{B.7}$$

For the VO1 EOS, convexity criteria (B.5) read,

$$\begin{cases} \text{(a)} & \rho P\left(\dfrac{R}{C_v}(1+a\rho) + \dfrac{1+2a\rho}{1+a\rho}\right) > 0, \\ \text{(b)} & \dfrac{P}{\rho RC_v(1+a\rho)} > 0, \\ \text{(c)} & -\dfrac{P}{C_v} < 0, \\ \text{(d)} & \dfrac{P^2(1+2a\rho)}{RC_v(1+a\rho)^2} > 0. \end{cases} \tag{B.8}$$

The VO1 EOS is consequently convex if,

$$a > -\frac{1}{\rho}. \tag{B.9}$$

In the present contribution, the virial coefficient a is taken positive to keep the VO1 EOS convex in any flow conditions.

## Appendix C. Parameters of the VO1 EOS with $C_v(T)$

This Appendix details the method for the determination of parameters $C_{v_0}$, c and q for the gas products resulting of reactive material combustion. The mixture internal energy reads,

$$e_{mix}(T) = Y\left(C_{v_0}T + \frac{c}{2}T^2 + q\right) + (1-Y)\left(C_{v_0,in}T + \frac{c_{in}}{2}T^2 + q_{in}\right). \tag{C.1}$$

For the sake of simplicity, a noble gas is used as inert. Its specific heat is constant on a wide range of temperature meaning that $c_{in} = 0$ and $C_{v_0,in} = C_{v,in}$ where $C_{v,in}$ is the constant specific heat known from the literature (Chase, 1998). Also, for noble gases their reference energy $q_{in}$ is zero. Equation (C.1) simplifies as,

$$e_{mix}(T) = Y\left(C_{v_0}T + \frac{c}{2}T^2 + q\right) + (1-Y)C_{v,in}T. \tag{C.2}$$

Multiple constant volume explosions are computed with CHEETAH for different mixtures i.e. different mass fraction Y. Equation (C.2) is expressed for each test run resulting in the following equation system, where subscript "j" denotes a specific mixture, $j = 1, \ldots, n$ where n is the number of tests considered.

$$C_{v_0}T_{flame,j} + \frac{c}{2}T_{flame,j}{}^2 + q = \frac{e_{mix,j}(T_{flame,j})}{Y_j} - \frac{(1-Y_j)}{Y_j}C_{v,in}T_{flame,j}. \tag{C.3}$$

At this point, only data computed with CHEETAH are used, meaning that heat losses are absent ($\pi = 0$). Consequently, for a constant volume explosion, energy conservation implies that the mixture energy is conserved throughout the whole reaction process,

$$e_{mix,j}(T_{flame,j}) = e_{mix,j}(T_i) = Y_j e_s^i + (1-Y_j)C_{v,in}T_0. \tag{C.4}$$

The solid reactive material initial energy $e_s^i$ is computed directly with CHEETAH for $Y = 1$, i.e. the reactive material alone. The initial temperature of the mixture is $T_0 = 298$ K. Inserting (C.4), the equation system (C.3) now reads,

$$C_{v_0}T_{flame,j} + \frac{c}{2}T_{flame,j}{}^2 + q = e_s^i - \frac{(1-Y_j)}{Y_j}C_{v,in}(T_{flame,j} - T_0). \tag{C.5}$$

$C_{v_0}$, c and q are deduced from system (C.5) with the least square method.

This method is only applicable when the composition of gas products resulting of reactive materials combustion is constant for every mixture used in the least square method. Otherwise, resulting parameters $C_{v_0}$ and c will likely be unphysical by implicitly trying to account for composition changes in the gaseous phase. A fine indicator of changes in combustion products is the molar mass $\widehat{W}$ of gases produced by the reactive material combustion. Considering various materials presented in Section 4 (NC-13, RDX, NG and HMX), Figure 9 shows molar mass $\widehat{W}$ of gas products resulting of reactive material combustion, when inert Argon is present, computed with the thermochemical code CHEETAH for mixtures between $Y = 15\%$ and $Y = 100\%$.

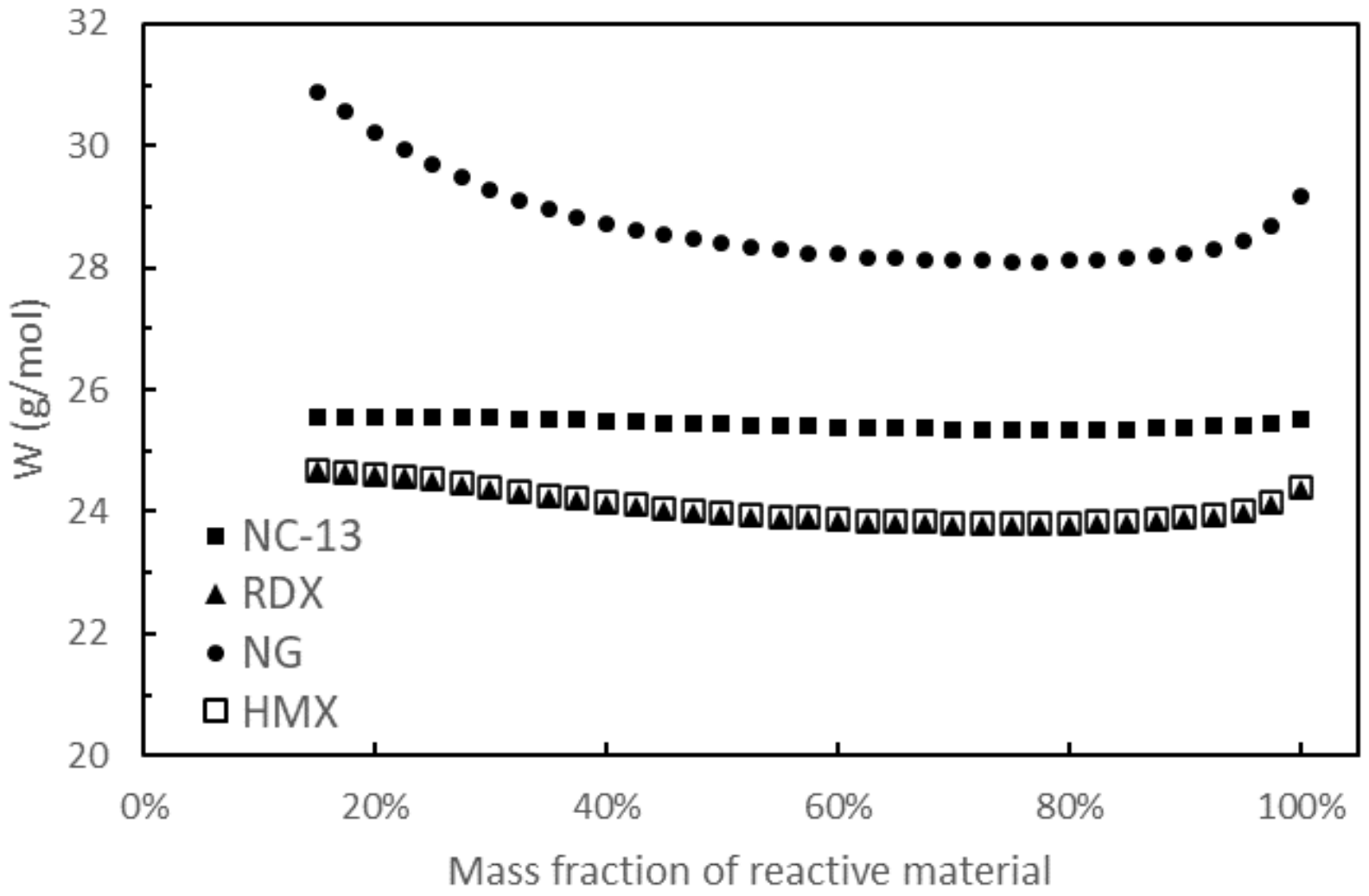

**Figure 9.** Molar mass of gas products resulting of reactive material combustion in inert Argon, computed with the thermochemical code CHEETAH for various mass fraction of reactive material. Considered materials are NC-13 in black filled squares, RDX in black filled triangles, NG in black filled circles and HMX in black empty squares. RDX and HMX results are merged.

Figure 9 shows that only gas products resulting of NC-13 combustion seems to be frozen when diluted in inert. Therefore, only NC-13 is retained among considered reactive materials in this paper to study the impact of a specific heat temperature dependence in the VO1 formulation.

Example of NC-13 variable specific heat parameters are given in Table III. For the determination of these parameters, 35 mixtures between $Y = 15\%$ and $Y = 100\%$ incrementing by 2.5% are considered, resulting in wide temperature variations. Mixtures below $Y = 15\%$ are not available with the thermochemical code CHEETAH as the pressure is too low leading to failed or poorly accurate computations. Flame temperature as a function of reactant mass fraction is shown in Figure 10. Each point of the graph corresponds to data indexed by j used in the least squares method to determine $C_{v_0}$, c and q with system (C.5).

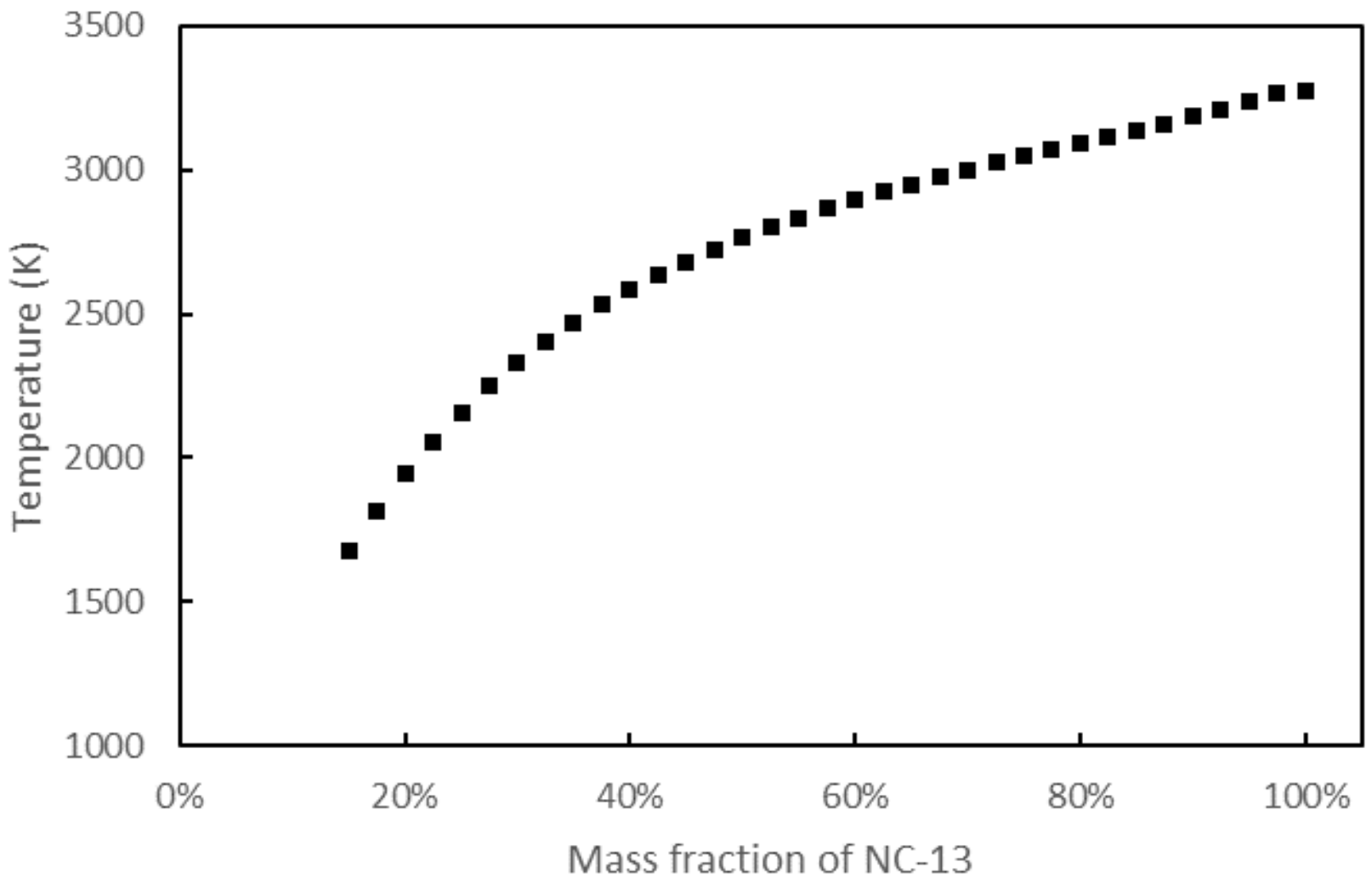

**Figure 10.** Flame temperature of various mixtures made of NC-13 and Argon. These mixtures are used in system (C.5) to determine parameters of the combustion products temperature dependent specific heat.

The temperature is ranging from approximately 3300 K, which is the flame temperature of pure NC-13, to 1600 K roughly half the flame temperature. This temperature range indicates validity domain of the specific heat parameters $C_{v_0}$ and c.

## Appendix D. Frozen mixture speed of sound calculation

Sound speed definitions (6.14) for MNA and (6.22) MVO1 EOS are derived in this Appendix. Considering the same ideal gas mixture evolving in temperature and pressure equilibrium of Section 6, the frozen sound speed is defined as,

$$c_{mix}^2 = \left( \frac{\partial P}{\partial \rho_{mix}} \right)_{s_{mix}, Y_1, \ldots, Y_N}, \tag{C.6}$$

where $s_{mix}$ represents the mixture specific entropy. The Gibbs identity adapted to a mixture is obtained as,

$$de_{mix} = Tds_{mix} - Pdv_{mix} + \sum_{k=1}^{N} \mu_k dn_k, \tag{C.7}$$

where $n_k$ and $\mu_k$ are respectively the mole number and the chemical potential of component k. However, a frozen mixture implies that $dn_k = 0$ for all $k = 1, \ldots, N$ components, which simplifies the Gibbs identity.

The differential of function $e_{mix}(P, \rho_{mix}, Y_1, \ldots, Y_N)$ for a constant composition reads,

$$de_{mix} = \left( \frac{\partial e_{mix}}{\partial P} \right)_{\rho_{mix}, Y_1, \ldots, Y_N} dP + \left( \frac{\partial e_{mix}}{\partial \rho_{mix}} \right)_{P, Y_1, \ldots, Y_N} d\rho_{mix}. \tag{C.8}$$

Combining expressions (C.7) and (C.8) and considering an isentropic process yields the sound speed definition,

$$c_{mix}^2 = \left( \frac{\partial P}{\partial \rho_{mix}} \right)_{s_{mix}, Y_1, \ldots, Y_N} = \frac{\dfrac{P}{\rho_{mix}^2} - \left( \dfrac{\partial e_{mix}}{\partial \rho_{mix}} \right)_{P, Y_1, \ldots, Y_N}}{\left( \dfrac{\partial e_{mix}}{\partial P} \right)_{\rho_{mix}, Y_1, \ldots, Y_N}}. \tag{C.9}$$

Following derivations focuses on demonstrating that definition (6.22) is correct. Definition (6.22) reads,

$$c_{mix}^2 = \frac{C_{p,mix}}{C_{v,mix}} \left( \frac{\partial P}{\partial \rho_{mix}} \right)_{T, Y_1, \ldots, Y_N}. \tag{C.10}$$

Heat capacities of the mixture are defined as,

$$C_{v,mix} = \left( \frac{\partial e_{mix}}{\partial T} \right)_{\rho_{mix}, Y_1, \ldots, Y_N} \quad \text{and} \quad C_{p,mix} = \left( \frac{\partial h_{mix}}{\partial T} \right)_{P, Y_1, \ldots, Y_N}. \tag{C.11}$$

Moreover, one can write,

$$\left( \frac{\partial P}{\partial \rho_{mix}} \right)_{T, Y_1, \ldots, Y_N} \left( \frac{\partial \rho_{mix}}{\partial T} \right)_{P, Y_1, \ldots, Y_N} \left( \frac{\partial T}{\partial P} \right)_{\rho_{mix}, Y_1, \ldots, Y_N} = -1. \tag{C.12}$$

Introducing expressions (C.11) and (C.12) in equation (C.10),

$$c^2_{mix} = -\frac{\left(\frac{\partial h_{mix}}{\partial \rho_{mix}}\right)_{P,Y_1,\ldots,Y_N}}{\left(\frac{\partial e_{mix}}{\partial P}\right)_{\rho_{mix},Y_1,\ldots,Y_N}}. \tag{C.13}$$

Knowing the specific enthalpy definition $h_{mix} = e_{mix} + Pv_{mix}$ (C.13) becomes,

$$c^2_{mix} = \frac{\frac{P}{\rho^2_{mix}} - \left(\frac{\partial e_{mix}}{\partial \rho_{mix}}\right)_{P,Y_1,\ldots,Y_N}}{\left(\frac{\partial e_{mix}}{\partial P}\right)_{\rho_{mix},Y_1,\ldots,Y_N}}, \tag{C.14}$$

which is exactly the sound speed definition (C.9).

## Author declarations

### Conflicts of interest

The authors have no conflicts to disclose.

## Data availability

The data that supports the findings of this study are available within the article.